\documentstyle[11pt,amsfonts,fullpage]{amsart}
 
\def\endproof{\qed\smallskip}
\def\blacksquare{\hbox to .60em{\vrule width .60em height .60em}}

\newtheorem{theorem}{Theorem}[section]
\newtheorem{corollary}[theorem]{Corollary}

\newtheorem{lemma}[theorem]{Lemma}
\newtheorem{proposition}[theorem]{Proposition}
\newtheorem{remark}[theorem]{Remark}

\begin{document}

\title[AH Einstein Metrics on 4-Manifolds]{Boundary Regularity, 
Uniqueness and Non-Uniqueness for AH Einstein Metrics on 4-Manifolds}

\author[M. Anderson]{Michael T. Anderson}

\thanks{Partially supported by NSF Grant DMS 0072591}

\maketitle

\abstract
This paper studies several aspects of asymptotically hyperbolic 
Einstein metrics, mostly on 4-manifolds. We prove boundary regularity 
(at infinity) for such metrics and establish uniqueness under natural 
conditions on the boundary data. By examination of explicit black hole 
metrics, it is shown that neither uniqueness nor finiteness holds in 
general for AH Einstein metrics with a prescribed conformal infinity. 
We then describe natural conditions which are sufficient to ensure 
finiteness.

\endabstract

\contentsline {section}{\tocsection {}{0}{Introduction.}}{1}
\contentsline {section}{\tocsection {}{1}{Conformally Compact Einstein Metrics.}}{4}
\contentsline {section}{\tocsection {}{2}{Boundary Regularity.}}{8}
\contentsline {section}{\tocsection {}{3}{Uniqueness.}}{13}
\contentsline {section}{\tocsection {}{4}{Non-Uniqueness.}}{19}
\contentsline {section}{\tocsection {}{5}{Cusp Formation and Hyperbolic Manifolds.}}{24}
\contentsline {section}{\tocsection {}{}{References}}{32}

\setcounter{section}{-1}

\section{Introduction.}
\setcounter{equation}{0}

 In this paper, we study several aspects of asymptotically hyperbolic 
(AH) Einstein metrics on an open 4-manifold $M$ with compact boundary 
$\partial M.$ These metrics are complete Einstein metrics $g$ on $M$, 
normalized so that
\begin{equation} \label{e0.1}
Ric_{g} = - 3g, 
\end{equation}
which are conformally compact in the sense of Penrose, in that there 
exists a defining function $\rho $ for $\partial M$ in $M$, such that 
the conformally equivalent metric
\begin{equation} \label{e0.2}
\bar g = \rho^{2}\cdot  g 
\end{equation}
extends to a Riemannian metric on $\bar M.$ Recall that a defining 
function $\rho $ for $\partial M$ is essentially just a coordinate 
function for $\partial M$ in $M$; thus $\rho $ is a smooth, typically 
$C^{\infty},$ function on $\bar M = M\cup\partial M$ such that $\rho  > 
$ 0 on $M$, $\rho^{-1}(0) = \partial M$ and $d\rho  \neq $ 0 on 
$\partial M.$ 

 Defining functions are unique only up to multiplication by positive 
functions on $\bar M.$ Hence only the conformal class $[\bar g]$ is 
uniquely determined by $g$, as is the conformal class $[\gamma ]$ of 
the induced metric $\gamma  = \bar g_{\partial M}$ on $\partial M.$ The 
class $[\gamma ]$ is called the {\it  conformal infinity} of $g$ and a 
choice $\gamma\in [\gamma ]$ will be called a {\it  boundary metric}. 
The metric $g$ is $L^{k,p}$ or $C^{m,\alpha}$ conformally compact if 
there exists a defining function such that the metric $\bar g$ in (0.2) 
has a $L^{k,p}$ or $C^{m,\alpha}$ extension to $\bar M;$ here $L^{k,p}$ 
is the Sobolev space of $k$ weak derivatives in $L^{p}$ and 
$C^{m,\alpha}$ is the usual H\"older space. It is easy to see, cf. \S 
1, that a conformally compact Einstein metric has curvature decaying to 
$- 1$ at an exponential rate, so that such manifolds are asymptotically 
hyperbolic.

 Regarding the existence of such metrics, Graham and Lee [19] have 
proved that any metric $\gamma $ near the standard metric $\gamma_{o}$ 
on $S^{n-1}$ in a sufficiently smooth topology may be filled in with an 
AH Einstein metric $g$ on the $n$-ball $B^{n}$ having prescribed 
boundary metric $\gamma .$ Further, such AH Einstein metrics have a 
conformal compactification with a certain degree of smoothness. More 
precisely, they prove that there is an open neighborhood ${\cal 
U}_{\gamma_{o}}$ of $\gamma_{o}$ in the space of $C^{m,\alpha}$ metrics 
${\cal M}^{m,\alpha}(S^{n-1})$ on $S^{n-1},$ for any $m \geq $ 2, such 
that any metric $\gamma\in{\cal U}_{\gamma_{o}}$ is the boundary metric 
of an AH Einstein metric $g$ on the $n$-ball $B^{n},$ i.e. $\bar 
g|_{\partial M} = \gamma .$ Further, the metric $g$ is $C^{n-2,\alpha}$ 
conformally compact for $n > $ 4 and $C^{1,\alpha}$ compact, for $n =$ 
4. 

 Recently, Biquard [7] has extended this result to boundary metrics 
$\gamma $ in an open neighborhood ${\cal U}_{\gamma_{o}} \subset  {\cal 
M}^{m,\alpha}(\partial M)$ of the boundary metric $\gamma_{o}$ of an 
arbitrary non-degenerate AH Einstein manifold $(M, g)$. Here 
non-degenerate means that there are no non-trivial $L^{2}$ 
infinitesimal AH Einstein deformations of $(M, g)$. Biquard's method 
can be shown to give a $C^{2}$ conformal compactification for $n \geq 
4$, (by choosing $\delta  =$ 2 in the notation of [7]).

\medskip

 The first purpose of this paper is to study the boundary regularity of 
conformal compactifications $\bar g$ of AHE metrics $g$, in dimension 
4. Namely, given an AH Einstein metric $(M, g)$ which, in some 
compactification $\bar g$ as in (0.2), has a $C^{m,\alpha}$ boundary 
metric $\gamma ,$ is there a $C^{m,\alpha}$ conformal compactification 
of $g$? This issue of boundary regularity was first raised in [16], and 
has been an open problem for some time. In fact, when $n =$ dim $M$ is 
odd, it was discovered in [16] that boundary regularity in general 
breaks down at the order $m = n- 1,$ in that there are log terms in the 
asymptotic expansion of $\bar g$ near $\partial M$ at this order, cf. 
(0.4) below. These log terms play an important role in the AdS/CFT 
correspondence relating string theory and conformal field theory, cf. 
[38], [25]. 

 When $n$ = dim $M$ is even, there are no log terms in the expansion, 
and it is possible that a $C^{\infty}$ boundary metric has a 
$C^{\infty}$ compactification. We resolve this issue in dimension 4.

\begin{theorem} \label{t 0.1.}
  Let $M$ be a 4-manifold and $g$ an AH Einstein metric on M. Suppose 
$g$ has an $L^{2,p}$ conformal compactification, for some $p > 4$, in 
which the boundary metric $\gamma $ is $C^{m,\alpha}$, $m \geq  3, 
\alpha > 0$. Then $(M, g)$ has a $C^{m,\alpha}$ conformal 
compactification with the same boundary metric. This result holds also 
if $m = \infty $ or $m = \omega .$
\end{theorem}

 The reason for considering the (weak) condition of $L^{2,p}$ 
compactification is that the result of Biquard above naturally gives 
the existence of AH Einstein metrics with such compactifications. 
Theorem 0.1 is proved in \S 2, cf. Theorem 2.4 and Corollary 2.5.

\medskip

 A second purpose of the paper is to study the uniqueness question, 
i.e. to what extent an AH Einstein metric $(M, g)$ is uniquely 
determined by its conformal infinity $(\partial M, [\gamma ]),$ or by 
data of $(M, g)$ at infinity. Again, this issue, raised in [16], has 
been open for some time, cf. also [8].

 We summarize some of the results on non-uniqueness here, and refer to 
\S 4 for detailed statements and constructions. First, it turns out 
that it has been known to physicists working in Euclidean quantum gravity 
for a rather long time that uniqueness fails in general. For example, there 
are distinct, i.e. non-isometric, AH Einstein metrics on $M = 
S^{2}\times {\Bbb R}^{2},$ which have the same conformal infinity on 
$\partial M = S^{2}\times S^{1}.$ These metrics are from the family of 
AdS-Schwarzschild metrics, and have been analysed in detail in a 
remarkable paper of Hawking and Page [23]. Further there are AH 
Einstein manifolds $(M, g)$ of distinct topological type, with the same 
conformal infinity, so that the conformal infinity $(\partial M, \gamma 
)$ does not determine the topological type of $M$.

 In fact, even in dimension 3, where Einstein metrics are hyperbolic, 
i.e. of constant curvature, there are numerous examples of 
non-uniqueness. These examples come from the Thurston theory of Dehn 
surgery, or the well-known process of "opening cusps", cf. [20], [37]. 
This construction gives infinitely many isometrically distinct 
hyperbolic 3-manifolds, all with a given conformal infinity. 

 Of course such Thurston-Dehn surgery is a special feature of 
hyperbolic 3-manifolds, and does not generalize to hyperbolic 
$n$-manifolds, $n \geq $ 4. However, we will see in \S 4 that such Dehn 
surgery constructions do generalize in the context of AH Einstein 
metrics. As in dimension 3, this gives rise, in special situations, to 
an infinite sequence distinct AH Einstein metrics on a fixed 
4-manifold, with a fixed conformal infinity, cf. Proposition 4.4 and 
Remark 4.5. These metrics again come from 'AdS black hole' metrics, 
with a toral $(T^{2})$ event horizon and have been extensively examined 
in the AdS/CFT correspondence. In particular, we see that even the 
finiteness issue, in addition to the uniqueness issue, fails in general.

 In sum, there are numerous counterexamples to uniqueness of AH 
Einstein metrics with a given conformal infinity. We give at least a 
brief overview of this situation in \S 4.

\medskip

 Thus, it is of interest to understand under what conditions one can 
uniquely characterize an AH Einstein metric. To do this, recall the 
Fefferman-Graham expansion of a conformal compactification. Thus, let 
$\bar g$ be a compactification as in (0.2) where the defining function 
$t$ has the property that $t(x) = dist_{\bar g}(x, \partial M),$ in a 
collar neighborhood $U$ of $\partial M.$ It is easy to see that, given 
a boundary metric $\gamma $ in the conformal infinity of $(M, g)$, 
there is a unique such defining function $t = t_{\gamma}$ having 
$\gamma $ as its boundary metric, (cf. \S 1). The Gauss Lemma then 
implies that the metric $\bar g$ splits in $U$ as
\begin{equation} \label{e0.3}
\bar g = dt^{2} + g_{t}, 
\end{equation}
where $g_{t}$ is a curve of metrics on $\partial M.$ It follows from 
Theorem 0.1, c.f Corollary 2.5 below, that if the boundary metric 
$\gamma \in C^{m+1,\alpha}$, then the curve $g_{t}$ is at least a 
$C^{m,\alpha}$ curve of metrics in $t$. Following [16], we may then 
expand $g_{t}$ as a Taylor series in $t$ as
\begin{equation} \label{e0.4}
g_{t} = g_{(0)} + tg_{(1)} + t^{2}g_{(2)} + t^{3}g_{(3)} + ... + 
t^{m}g_{(m)} + O(t^{m+\alpha}). 
\end{equation}
One has $g_{(0)} = \gamma , g_{(1)} = 0$, while $g_{(2)}$ is determined 
locally from the geometry of the boundary metric $\gamma .$ On the 
other hand, the terms $g_{(j)}$ for $j \geq 3$ are not locally 
determined by the boundary $\gamma $ in general, cf. \S 3 for further 
discussion.

\begin{theorem} \label{t 0.2.}
  Suppose dim $M =$ 4 and the boundary metric $\gamma\in C^{7,\alpha}.$ 
Then the data $(\gamma , g_{(3)})$ on $\partial M$ uniquely determine 
an AH Einstein metric up to local isometry, i.e. if $g^{1}$ and $g^{2}$ 
are two AH Einstein metrics on manifolds $M^{1}$ and $M^{2},$ with 
$\partial M^{1} = \partial M^{2} = \partial M$ such that, w.r.t. 
geodesic compactifications as in (0.3),
\begin{equation} \label{e0.5}
\gamma^{1} = \gamma^{2} \ \ {\rm and} \ \ g_{(3)}^{1} = g_{(3)}^{2}, 
\end{equation}
then $g^{1}$ and $g^{2}$ are locally isometric and the manifolds 
$M^{1}$ and $M^{2}$ are commensurable, i.e. they have diffeomorphic 
universal covers.
\end{theorem}

 It follows that $M^{1}$ is diffeomorphic to $M^{2}$ and $g^{1}$ is 
isometric to $g^{2}$ if $\pi_{1}(M^{1}) \cong  \pi_{1}(M^{2})$ and the 
actions of $\pi_{1}(M^{i})$ on the universal cover $\widetilde M$ of 
$M^{i}$ are conjugate in the isometry group of $\widetilde M.$ Theorem 
0.2 is proved in \S 3.

\medskip

 The third purpose of the paper is to analyse the finiteness issue, 
i.e. when a given boundary metric $(\partial M, \gamma )$ bounds 
finitely or infinitely many AH Einstein 4-manifolds $(M_{i}, g_{i}).$ 
The main result is, roughly speaking, that under reasonably natural 
conditions only conformally flat boundary data can bound infinitely 
many AH Einstein metrics. Since the exact statement is somewhat 
technical, we refer to Theorem 5.3 for details. We also mention here 
Proposition 5.1, which gives a very simple proof of the result of 
Witten-Yau [39] on the connectedness of $\partial M$ when $\partial M$ 
has a component of positive scalar curvature.

\medskip

 This paper does not address the existence question for AH Einstein 
metrics, i.e. given a conformal class $[\gamma ]$ of metrics on 
$\partial M,$ does there exist an AH Einstein metric on $M$ which has 
$[\gamma ]$ as its conformal infinity. This issue will be discussed in 
a sequel [4] to this paper; the results obtained here however will play 
an important role in the sequel.

 As indicated above, there is a very extensive and active recent 
physics literature on AH Einstein metrics in relation to the AdS/CFT 
correspondence and this work is a strong influence on this paper. We 
refer to [14], [24], [31], [33], [38] for some relevant perspectives.

  I'd like to thank Robin Graham and Rafe Mazzeo for their comments and for 
pointing out a gap in an earlier proof of Proposition 2.1. Many thanks also 
to the referees for their suggestions in improving the exposition of the 
paper.

\section{Conformally Compact Einstein Metrics.}
\setcounter{equation}{0}

 In this section, we discuss some background material for conformally 
compact Einstein metrics. Although the results of this section hold in 
arbitrary dimensions, we carry out the computations only in dimension 
4, since this is the most relevant case for the paper; cf. also Remark 
1.5 for the higher dimensional case. Thus, we assume that $M$ is an 
open, connected, oriented 4-manifold, with non-empty and compact 
boundary. It is not assumed that $\partial M$ is connected.

 Let $\rho $ be a defining function for $\partial M$ and as in (0.2), 
set
\begin{equation} \label{e1.1}
\bar g = \rho^{2}\cdot  g. 
\end{equation}
Unless $\rho $ is restricted by the geometry of $g$, without loss of 
generality we may, and do, assume that $\rho $ is $C^{\infty}$ on $\bar 
M,$ so that the metric $\bar g$ is $C^{\infty}$ in the interior of $M$. 
The metric $\bar g$ is a $L^{k,p}$ or $C^{m,\alpha}$ compactification 
if it extends to a $L^{k,p}$ or $C^{m,\alpha}$ metric on $\bar M.$ 
Thus, there are coordinate charts for a neighborhood of $\partial M$ in 
$M$ such that the local components of $\bar g$ in these charts are 
$L^{k,p}$ or $C^{m,\alpha}$ functions of the coordinates. We will 
usually assume that 
\begin{equation} \label{e1.2}
k \geq  2, \ p >  4, \ \ {\rm or} \ \ m+\alpha  >  1. 
\end{equation}
Sobolev embedding, in dimension 4, implies that for $p > $ 4, $L^{2,p} 
\subset  C^{1,\alpha}, \alpha  = 1-\frac{4}{p} > $ 0, while for $p > $ 
2, $L^{2,p} \subset  C^{\alpha}, \alpha  = 2-\frac{4}{p}.$ It is 
well-known that metrics have optimal regularity in harmonic 
coordinates. Such coordinates exist on $(\bar M, \bar g)$ if $\bar g$ 
is a $L^{k,p}$ metric with $k \geq $ 2, $p > $ 2; further the 
components of $\bar g$ are in $L^{k,p}$ w.r.t. such harmonic 
coordinates, cf. [6, Ch.5E].

 If $\bar g$ is a $C^{m,\alpha}$ compactification, then the boundary 
metric $\gamma  = \bar g|_{\partial M}$ on $\partial M$ is also 
$C^{m,\alpha},$ while if $\bar g$ is $L^{k,p},$ then the boundary 
metric $\gamma $ is $L^{k-\frac{1}{p}, p}$, cf. [1, 7.56]. In particular, 
by Sobolev embedding, a $L^{2,p}$ compactification has a $C^{1,\alpha}$ 
boundary metric when $p > 4$. Of course the converse of these statements may 
not hold in general. The degree of smoothness of the boundary metric 
$\gamma $ does not, apriori, imply any degree of smoothness of the 
compactification.

 As noted in \S 0, defining functions $\rho $ on $M$ are not unique, 
but differ by multiplication by positive functions. Conversely, given 
any positive smooth function $\phi $ on $M$ and any defining function 
$\rho ,$ then $\phi\cdot \rho $ is a defining function. Hence only the 
conformal class of the boundary metric $\gamma  = \bar g_{\partial M}$ 
is uniquely determined by $(M, g)$. 

 The curvatures of the Einstein metric $g$ and a compactification $\bar 
g$ are related by the following formulas in dimension 4:
\begin{equation} \label{e1.3}
\bar K_{ab} = \frac{K_{ab} + |\bar \nabla\rho|^{2}}{\rho^{2}} -  
\frac{1}{\rho}\{\bar D^{2}\rho (\bar e_{a},\bar e_{a})+\bar D^{2}\rho 
(\bar e_{b},\bar e_{b})\}. 
\end{equation}
\begin{equation} \label{e1.4}
\bar Ric = - 2\frac{\bar D^{2}\rho}{\rho} + (3\rho^{-2}(|\bar 
\nabla\rho|^{2}- 1) -  \frac{\bar \Delta\rho}{\rho})\bar g, 
\end{equation}
\begin{equation} \label{e1.5}
\bar s = - 6\frac{\bar \Delta\rho}{\rho} + 12\rho^{-2}(|\bar 
\nabla\rho|^{2} -  1). 
\end{equation}
The terminology here is the following: $D^{2}$ is the Hessian, $\nabla$ is 
the gradient, $\Delta = tr D^{2}$ the Laplacian, while $K_{ab}$ denotes 
sectional curvature and $\{e_{a}\}$ an orthonormal basis. (See [6, Ch.1J] 
for example for formulas for conformal changes of metric). All barred 
quantities are w.r.t. the $\bar g$ metric. The equation (1.4) is equivalent 
to the Einstein equations (0.1). Similar formulas hold in all dimensions. 
We also let
\begin{equation} \label{e1.6}
r = \log (\frac{2}{\rho}), \ \ \rho  = 2e^{- r.} 
\end{equation}

 Since $\rho $ is smooth on $\bar M,$ it is essentially immediate from 
(1.4) that if $\bar g$ is a $L^{k,p}$ compactification, satisfying 
(1.2), then $|\bar \nabla\rho| \equiv $ 1 at $\partial M,$ and hence by 
(1.3), the sectional curvatures of $g$ tend to $- 1$ at infinity in 
$(M, g)$. Hence any such conformally compact Einstein manifold is 
asymptotically hyperbolic (AH). Further at $\partial M,$ we have $\bar 
D^{2}\rho  = A$, where $A$ is the $2^{\rm nd}$ fundamental form of 
$\partial M$ in $(\bar M, \bar g);$ $A$ is $C^{\alpha}$ on $\partial 
M,$ again by (1.2). The equation (1.4) further implies that $\partial 
M$ is {\it  umbilic}, i.e. $A = \lambda\cdot \gamma ,$ for some 
function $\lambda $ on $\partial M.$ 

 From the formulas (1.4)-(1.5), it is clear that defining functions 
which satisfy
\begin{equation} \label{e1.7}
|\bar \nabla\rho| \equiv  1 
\end{equation}
in a collar neighborhood $U$ of $\partial M$ in $M$ are especially 
natural. Such defining functions will be called {\it geodesic}  
defining functions, (although they are called special defining 
functions in [18]). A brief computation shows that, (in general), 
\begin{equation} \label{e1.8}
|\bar \nabla\rho| = |\nabla r|, 
\end{equation}
where the norm on the left is w.r.t. $\bar g$ and that on the right 
w.r.t. $g$. Thus the condition (1.7) is an intrinsic property of $(M, 
g)$ and $r$. The function $r$ is a (signed) distance function on $(M, 
g)$ and the integral curves of $\nabla r$ are geodesics in $(M, g)$. 
Similarly $\rho $ is the distance function from $\partial M$ w.r.t. 
$\bar g.$ Geodesic defining functions are geometric, in that they 
depend on the geometry of $(M, \bar g)$ or $(M, g)$, and so their 
smoothness depends on the metric $\bar g.$ Thus, such functions will 
not be $C^{\infty}$ unless the compactification $\bar g$ is; if $\bar 
g$ is $C^{m,\alpha},$ then $\rho $ is $C^{m+1,\alpha}$ off the cutlocus 
of $\partial M$ in $(\bar M, \bar g).$

 Suppose there is a compactification $\widetilde g = \rho^{2}g$ of $(M, 
g)$ which is at least $C^{2},$ (actually $C^{1,1}$ suffices), with 
boundary metric $\gamma .$ Then it is easy to see that there is a 
unique geodesic defining function $t = t(\gamma )$ for $(M, g)$ such 
that the compactification
\begin{equation} \label{e1.9}
\bar g = t^{2}\cdot  g 
\end{equation}
has boundary metric $\gamma ,$ cf. [19, Lemma 5.2]. Briefly, write $t 
= u\cdot \rho $ where $u$ is a positive function on $\bar M$ with $u 
\equiv $ 1 on $\partial M$ so that the boundary metric of (1.9) is 
indeed $\gamma .$ Then the equation that $|\bar \nabla t|_{\bar g} =$ 
1, i.e. $t$ is a geodesic defining function, is equivalent to
\begin{equation} \label{e1.10}
2(\widetilde \nabla \rho )(\log u) + \rho|\widetilde \nabla \log 
u|_{\widetilde g}^{2} = \rho^{-1}(1 -  |\widetilde \nabla 
\rho|_{\widetilde g}^{2}). 
\end{equation}
This is a non-characteristic $1^{\rm st}$ order PDE, with $C^{2}$ (or 
$C^{1,1})$ coefficients, with right hand side in $C^{1}$ (or 
Lipschitz). Hence the Cauchy problem has a unique solution with $u 
\equiv  1$ on $\partial M$ in a collar neighborhood $U$ of $\partial 
M$. Observe however that $\bar g$ may not be as smooth as $\widetilde 
g$; if $\widetilde g \in C^{m,\alpha}$, $\alpha \geq 0$, then we have 
only $u \in C^{m-1,\alpha}$ and so $\bar g \in C^{m-1,\alpha}$ 

  The Gauss Lemma implies that the metrics $\bar g$ and $g$ split in 
$U$, as in (0.3):
\begin{equation} \label{e1.11}
\bar g = dt^{2} + g_{t}, \ \ {\rm and} \ \   g = dr^{2} + g_{r}, 
\end{equation}
with $g_{o} = \gamma $ and $g_{r} = t^{-2}g_{t}.$ The 1-parameter 
family $g_{t}$ is a $C^{m,\alpha}$ smooth curve of metrics on $\partial 
M$ if $\bar g$ is a $C^{m,\alpha}$ compactification. Observe also that 
the $2^{\rm nd}$ fundamental form $\bar A$ of the level sets of $t$ is 
given by $\bar A = \bar D^{2}t$ and so in particular by (1.4)
\begin{equation} \label{e1.12}
\bar A = 0 \ \ {\rm at} \ \ \partial M, 
\end{equation}
i.e. $\partial M$ is totally geodesic in $M$. In sum, if $\widetilde g$ 
is a $C^{m,\alpha}$ compactification with $m \geq $ 2 and $\alpha  \geq 
$ 0, then there is a unique compactification $\bar g$, at least 
$C^{m-1,\alpha}$, by a geodesic defining function inducing the boundary 
metric $\gamma $ of $\widetilde g$ on $\partial M.$ Such a 
compactification will be called the {\it  geodesic compactification}  
associated with $\gamma .$

\medskip

 If $\widetilde g$ is only $L^{2,p}$ or $C^{1,\alpha},$ then the 
discussion above does not hold; although geodesics normal to $\partial 
M$ do exist, they are not necessarily uniquely defined, and so one does 
not obtain a splitting (1.11) valid up to the boundary $\partial M.$ 
Nevertheless, the discussion above does hold "to first order" at 
$\partial M,$ in that there exists (another) $L^{2,p}$ compactification 
$\hat g$ which satisfies (1.12). 

\begin{lemma} \label{l 1.1.}
  Let $\widetilde g$ be a $L^{2,p}$ compactification of g, $p > 4$, 
with boundary metric $\gamma$ and $L^{3,p}$ defining function $\rho$. 
Then there exists another (possibly equal) $L^{2,p}$ compactification 
$\hat g$ of g, with the same boundary metric $\gamma$, such that
$$\hat A = 0 \ \ {\rm at} \ \ \partial M. $$
\end{lemma}

\noindent
{\bf Proof:}
 Let $\phi $ be a $L^{2,p}$ positive function on $\bar M,$ with $\phi  
\equiv $ 1 on $\partial M$ and set $\hat g = \phi^{2}\widetilde g.$ 
Then $\hat g \in  L^{2,p},$ the boundary metric of $\hat g$ is $\gamma 
,$ and the $2^{\rm nd}$ fundamental forms $\hat A$ and $\widetilde A$ 
of $\partial M$ w.r.t $\hat g$ and $\widetilde g$ are related by
$$\hat A = \widetilde A + <\widetilde \nabla \log\phi , \widetilde 
\nabla \rho>\cdot \gamma . $$
Let $\widetilde A$ be the $2^{\rm nd}$ fundamental form of the $\rho$-level 
sets w.r.t. $\widetilde g$ and set $\lambda = 
tr_{\widetilde g} \widetilde A / 3$, so that $\lambda \in L^{1,p}$. 
As noted above, $\widetilde A = \lambda\gamma$ and 
$\widetilde \nabla \rho = N$, the $\widetilde g$ unit normal, both at 
$\partial M$. Choosing $\phi $ to satisfy 
$<\widetilde \nabla \log\phi , \widetilde \nabla \rho> = 
N(\log \phi) = -\lambda$ at $\partial M$ then gives the result.

 Observe that if $\widetilde g\in L^{k,p}$ or $C^{m,\alpha}$ then $\phi$ 
may also be chosen to be in $L^{k,p}$ or $C^{m,\alpha}$.

{\endproof}

 The next result shows that the Ricci curvature $\bar Ric$ at $\partial 
M$ of a $C^{2}$ geodesic compactification $\bar g$ is determined by the 
intrinsic $C^{2}$ geometry of the boundary metric $\gamma .$

\begin{lemma} \label{l 1.2.}
  Let $\bar g$ be a $C^{2}$ geodesic compactification of an AH Einstein 
4-manifold $(M, g)$, with $C^{2}$ boundary metric $\gamma .$ Then at 
$\partial M,$
\begin{equation} \label{e1.13}
\bar s = 6\bar Ric(N, N) = \frac{3}{2}s_{\gamma}, 
\end{equation}
where $N$ is the unit normal to $\partial M$ w.r.t. $\bar g.$ If $X$ is 
tangent to $\partial M,$ then
\begin{equation} \label{e1.14}
\bar Ric(N, X) = 0, 
\end{equation}
while if $T$ denotes the projection onto $T(\partial M),$ then
\begin{equation} \label{e1.15}
(\bar Ric)^{T}= 2Ric_{\gamma}-  \frac{1}{4}s_{\gamma}\cdot \gamma . 
\end{equation}
\end{lemma}

\noindent
{\bf Proof:}
 The equality (1.14) follows immediately from (1.4), while the first 
equality in (1.13) follows from (1.4) and (1.5). For the rest, let $t$ 
be the geodesic defining function, let $\bar e_{i}, i =$ 1,2,3, be an 
o.n. basis at any point for a level set $S(t)$ of $t$, for $t$ near 0, 
and let $N = \bar \nabla t$ be the unit normal field to $S(t)$. By 
definition, the curvature $\bar K_{Ni} = \bar R(\bar e_{i}, N, N, \bar e_{i}) 
= <\bar \nabla_{\bar e_{i}}\bar \nabla_{N} N, \bar e_{i}> - 
<\bar \nabla_{N}\bar \nabla_{\bar e_{i}} N, \bar e_{i}> - 
<\bar \nabla_{[\bar e_{i}, N]} N, \bar e_{i}>$. Using the fact that $N$ is
tangent to geodesics and $N$ is a gradient, together with the fact that 
$\bar D^{2}t = A = 0$ at $\partial M$ by (1.12), gives
$$\bar K_{Ni} = - 
\lim_{t \rightarrow 0}t^{-1}\bar D^{2}t(\bar e_{i}, \bar e_{i}),$$
at $\partial M$. From (1.4)-(1.5), it then follows that 
\begin{equation} \label{e1.16}
\bar K_{Ni} =  \frac{1}{2}\bar Ric(\bar e_{i}, \bar e_{i}) -  
\frac{1}{12}\bar s, 
\end{equation}
again on $\partial M$. On the other hand, the Gauss-Codazzi equations and 
(1.12) again imply
$$\bar Ric(\bar e_{i}, \bar e_{i}) = Ric_{\gamma}(\bar e_{i}, \bar 
e_{i}) + \bar K_{Ni}, $$
and so substituting in (1.16) gives $\frac{1}{2}\bar Ric(\bar e_{i}, 
\bar e_{i}) = Ric_{\gamma}(\bar e_{i}, \bar e_{i}) -  \frac{1}{12}\bar 
s.$ Thus (1.15) follows from the second equality in (1.13).

 Finally, to prove this second equality, Gauss-Codazzi and (1.12) again give 
$\frac{1}{2}s_{\gamma} = \sum_{i<j\leq 3}\bar K_{ij}$. But by definition,
$$\frac{1}{2}\bar s = \sum_{i\leq 3}\bar K_{Ni}+\sum_{i<j\leq 3}\bar 
K_{ij}.$$
The first term is just $\bar Ric(N,N)$, and so $\frac{1}{2}\bar s = 
\bar Ric(N,N) + \frac{1}{2}s_{\gamma}$. From (1.4) and (1.5), this gives 
$\frac{1}{2}\bar s = \frac{1}{6}\bar s + \frac{1}{2}s_{\gamma}$, 
which gives (1.13).
{\endproof}

  We will also need the analogue of Lemma 1.2 when the compactification is 
not geodesic.

\begin{lemma} \label{l 1.3.}
Let $\bar g$ be a geodesic compactification and let $\widetilde g = 
\phi^{2}\bar g$ be another compactification with the same boundary metric 
$\gamma$, so that $\phi \equiv 1$ on $\partial M$. Suppose that 
$\widetilde A = 0$ on $\partial M$, cf. Lemma 1.1. Then the Ricci curvature 
$\widetilde Ric$ of $\widetilde g$ at $\partial M$ is determined by the 
Ricci curvature of the boundary metric $\gamma$ and the scalar curvature 
$\widetilde s$ of $\widetilde g$ at $\partial M$. In fact, at $\partial M$,
\begin{equation} \label{e1.17}
\widetilde Ric = \bar Ric + \tfrac{1}{6}(\widetilde s - \bar s)(\widetilde g + 2\nu\cdot \nu),
\end{equation}
where $\bar Ric$ and $\bar s$ are as in (1.13)-(1.15) and $\nu$ is the unit 
conormal at $\partial M$. If $\widetilde A \neq 0$ at $\partial M$, then (1.17) holds modulo terms of the form $\widetilde A^{2}$.
\end{lemma}

\noindent
{\bf Proof:}
By standard formulas for conformal changes of the metric, cf. [6, 1.159], the 
Ricci curvatures of $\widetilde g$ and $\bar g$ are related by
\begin{equation} \label{e1.18}
\widetilde Ric = \bar Ric - 2\phi^{-1}\bar D^{2}\phi + 4(d\log\phi)^{2} - 
\phi^{-1}\bar \Delta \phi \cdot \bar g - |d\log \phi|^{2}\cdot \bar g.
\end{equation}
At $\partial M$, $\phi = 1$ and $d\phi = 0$, since $\widetilde A = 0$. Hence 
$$\widetilde Ric = \bar Ric - 2\bar D^{2}\phi - \bar \Delta \phi \cdot \bar g.$$ 
In an orthonormal frame $e_{1}, ..., e_{4}$, $e_{4} = N$, 
$\bar D^{2}\phi(e_{i}, e_{j}) = 0$ at $\partial M$ unless $e_{i} = e_{j} = N$,
 i.e. $\bar D^{2}\phi = [NN(\phi)]\nu \cdot \nu$. Hence $\bar D^{2}\phi(N,N) = 
\bar \Delta \phi$, and it follows that 
$$\widetilde Ric = \bar Ric - (\bar\Delta\phi)(\widetilde g + 2\nu\cdot \nu),$$ at $\partial M$. Taking the trace of this equation gives 
$\widetilde s = \bar s - 6\bar \Delta \phi$, which implies (1.17). If 
$\widetilde A \neq 0$, then the same arguments show that (1.18) gives (1.17) 
modulo $\widetilde A^{2}$ terms.
{\endproof}

 Next we have an interesting estimate for the scalar curvature $\bar s$ 
of geodesic compactifications $\bar g.$

\begin{proposition} \label{p 1.4.}
  Let $\bar g$ be a $C^{2}$ geodesic compactification with boundary 
metric $\gamma $ and scalar curvature $\bar s.$ Then off the cut locus 
of $t$ in (M, $\bar g),$
\begin{equation} \label{e1.19}
\bar s'  \equiv  <\bar \nabla \bar s, \bar \nabla t>  = 6t^{-1}|\bar 
D^{2}t|^{2} \geq  0. 
\end{equation}
In particular, $\bar s$ is uniformly bounded below in this region by 
its boundary value on $\partial M,$ and thus bounded below by the 
$C^{2}$ geometry of $\gamma .$
\end{proposition}

\noindent
{\bf Proof:}
 The computations below are w.r.t. the $\bar g$ metric, but we drop the 
bar from the notation. The flow lines of $\nabla t$ are geodesics, and 
hence a standard result in Riemannian geometry, (cf. [35] for example) 
implies that the following Ricatti equation holds in $U$:
\begin{equation} \label{e1.20}
H'  + |A|^{2} + Ric(\nabla t, \nabla t) = 0; 
\end{equation}
here $H$ = tr $A$ = $\Delta t$ and $H'  = <\nabla H, \nabla t> .$ By 
(1.4), we have $Ric(\nabla t, \nabla t) = - 2t^{-1}(D^{2}t)(\nabla t, 
\nabla t) -  t^{-1}\Delta t,$ and since $|\nabla t| =$ 1, 
$(D^{2}t)(\nabla t, \nabla t) =$ 0. Hence dividing (1.20) by $t$ gives
\begin{equation} \label{e1.21}
t^{-1}(\Delta t)'  + t^{-1}|D^{2}t|^{2} -  t^{-2}\Delta t = 0. 
\end{equation}
But $t^{-1}(\Delta t)'  = (t^{-1}\Delta t)'  + t^{-2}\Delta t,$ and so 
(1.21) becomes
$$(t^{-1}\Delta t)'  + t^{-1}|D^{2}t|^{2} = 0. $$
The equation (1.19) then follows from (1.5).
{\endproof}
\begin{remark} \label{r 1.5.}
  {\rm All of the discussion in this section holds in any dimension, 
with obvious modifications for Sobolev embedding and constants 
depending on dimension. For instance, in dimension $n$,  with 3 in 
(0.1) replaced by $n-1$, (1.13) holds with $\frac{3n-4}{2(n-1)(n-2)}$ 
in place of $\frac{3}{2},$ while the factor of 6 in (1.13) and (1.19) 
should be replaced by $2(n-1)$.}
\end{remark}

\section{Boundary Regularity.}
\setcounter{equation}{0}

 In this section, we study the boundary regularity of AH Einstein 
metrics on 4-manifolds and establish Theorem 0.1.

 As is well-known [6, Ch.5], Einstein metrics satisfy an elliptic 
system of equations in harmonic coordinates, and so one obtains higher 
order $(C^{\infty}$ or $C^{\omega})$ regularity of such metrics from a 
local $L^{p}$ bound on the curvature. With regard to boundary 
regularity, the boundary of an AH Einstein metric occurs at infinity. 
If one works in local coordinates for $\partial M,$ the system of 
Einstein equations becomes degenerate at $\partial M,$ and thus 
difficult to deal with for regularity issues.

 It is a special feature of dimension 4 that Einstein metrics $(M, g)$ 
also satisfy a conformally invariant equation, namely the Bach equation
\begin{equation} \label{e2.1}
\delta d(Ric -  \frac{s}{6}g) +\stackrel {\circ}W(Ric) = 0, 
\end{equation}
cf.  [6, (4.77)]. This is the Euler-Lagrange equation for ${\cal W} ,$ 
the square of the $L^{2}$ norm of the Weyl curvature $W$. Here $Ric$ 
and $g$ are viewed as 1-forms with values in $TM$, $d = d^{\nabla}$ is 
the exterior derivative $\Lambda^{1} \rightarrow  \Lambda^{2}$ defined 
in terms of the metric $g$ and $\stackrel {\circ} W$ is the action of 
the Weyl tensor on symmetric bilinear forms.

 Of course Einstein metrics satisfy (2.1) in any dimension, but the 
expression (2.1) is conformally invariant only in dimension 4. Hence 
(2.1) holds for any conformal compactification $(M, \bar g)$ of $(M, 
g)$. Observe that (2.1) is a $4^{\rm th}$ order system of equations in 
the metric $g$, as opposed to the $2^{\rm nd}$ order system of Einstein 
equations. Note also that, being conformally invariant, the equation 
(2.1) is trace-free, i.e. its trace vanishes identically.

 We first prove boundary regularity of an $L^{2,p}$ compactification, 
$p > $ 4, given suitable control on the scalar curvature of the 
compactification. As mentioned in the Introduction, we consider 
compactifications which are $L^{2,p},$ since this level of regularity 
exists for the AH Einstein metrics constructed by Biquard [7].

\begin{proposition} \label{p 2.1.}
  Let $(M, g)$ be an AH Einstein 4-manifold which admits an $L^{2,p}$ 
conformal compactification $\bar g,$ with boundary metric $\gamma ,$ 
for some $p > 4$, so that
\begin{equation} \label{e2.2}
|R_{\bar g}|_{L^{p}} \leq  \Lambda  <  \infty , 
\end{equation}
where $R$ is the curvature tensor.

 Let $k \geq 1$ and $q \geq 2$. If $\gamma  \in  L^{k+2,q}(\partial 
M)$, the scalar curvature $\bar s \in  L^{k,q}(\bar M)$, with 
$\bar s |_{\partial M} \in L^{k,q}(\partial M)$, then the 
metric $\bar g$ is in $L^{k+2,q}(\bar M)$. (This last assumption may be 
realized by assuming $\bar s \in L^{k+\frac{1}{q}, q}$, c.f. [1, 7.56]).

 The same result holds with respect to the H\"older $C^{m,\alpha}$ spaces, 
i.e. if $\gamma \in C^{m+2, \alpha}(\partial M)$ and 
$\bar s \in C^{m, \alpha}(\bar M)$, then $\bar g \in C^{m+2, \alpha}(\bar M)$. 
If $\bar s$ and $\gamma $ are $C^{\omega},$ i.e. real-analytic, then so is 
$\bar g.$

 Further, if, with respect to a fixed harmonic coordinate system for $\bar M$, 
\begin{equation} \label{e2.3}
||\gamma||_{L^{k+2,q}(\partial M)} \leq  C \ \ {\rm and} \ \  ||\bar 
s||_{L^{k,q}(\bar M)} + ||\bar s||_{L^{k,q}(\partial M)} \leq  C, 
\end{equation}
then 
\begin{equation} \label{e2.4}
||\bar g||_{L^{k+2,q}(\bar M)}\leq  C_{1}, 
\end{equation}
where $C_{1}$ depends only on C, $\Lambda ,$ an upper diameter bound 
and a lower bound for the geodesic ball volume ratio $vol_{\bar 
g}B_{x}(s)/s^{4}, x\in M, s \leq $ 1, on (M, $\bar g).$ The analogous 
estimate holds with respect to the $C^{m,\alpha}$ H\"older norms.
\end{proposition}

\noindent
{\bf Proof:}
 The idea of the proof is to apply boundary regularity results for 
elliptic systems, in connection with the Bach system (2.1). However, 
the leading operator $\delta d$ in (2.1) is not elliptic. This can be 
rectified by considering the operator $\delta d + 2\delta^{*}\delta .$ 
In fact a standard Weitzenbock formula gives
$$\delta d + 2\delta^{*}\delta  = 2D^{*}D + {\cal R} , $$
where $D^{*}D$ is the rough Laplacian, and ${\cal R} $ is a curvature 
term, cf. [5, p.288] or [6, (4.71)]. (The exact form of ${\cal R} $ is of no 
importance here). The Bianchi identity $\delta Ric = -\frac{1}{2}ds$, gives 
$2\delta^{*} \delta Ric = - D^{2}s$, while a straightforward computation, 
cf. again [5], shows $\delta d(sg) = -2 \Delta s \cdot g + 2D^{2}s$. Hence 
(2.1) can 
be rewritten as
\begin{equation} \label{e2.5}
2D^{*}DRic = - \frac{2}{3}D^{2}s - \frac{1}{3}\Delta s\cdot g + {\cal 
R}_{1}, 
\end{equation}
where ${\cal R}_{1} = -{\cal R} (Ric)- \stackrel {\circ}W(Ric)$ is a 
term quadratic in curvature. Here and in the following, all metric 
quantities are w.r.t. $\bar g$ but the overbar is omitted from notation. 

 The assumptions on $s$ and $\gamma$ give control on the right side of 
(2.5), while the left side of (2.5) is essentially an elliptic $4^{\rm th}$ 
order system in the metric $g$. In principle, the result then follows from 
general elliptic boundary regularity theory, but there are a fair number of 
details to address. To begin, as noted in \S 1, one may choose local harmonic 
coordinates for $\bar M$ in which the components of $g$, (i.e. $\bar g$), 
are $L^{2,p}.$ With respect to such coordinates, the components of the 
Ricci curvature are
\begin{equation} \label{e2.6}
Ric_{ab} = \frac{1}{2}\Delta g_{ab} + [Q_{1}(g, \partial g)]_{ab}, 
\end{equation}
where $Q_{1}$ is of lower order; $Q_{1}$ involves only quadratic 
expressions in $g$, $g^{-1}$ and $\partial g.$ Similarly, the rough 
Laplacian $D^{*}D$ is also the function Laplacian to leading order, in 
that, in harmonic coordinates
\begin{equation} \label{e2.7}
- (D^{*}DRic)_{ab} = \Delta (Ric_{ab}) + [Q_{3}(g, \partial^{j}g)]_{ab}, 
\end{equation}
where $Q_{3}$ involves only derivatives of $g$ up to order 3; we refer 
to [2, p.234] for instance for the exact calculations. Hence, in local 
harmonic coordinates, the left side of (2.5) has the form of the 
bi-Laplacian of the metric, $\Delta\Delta (g_{ab}),$ to leading order. 
Further, in such coordinates, the Laplacian has the form 
\begin{equation} \label{e2.8}
\Delta  = g^{kl}\partial_{k}\partial_{l}, 
\end{equation}
and so involves the metric only to $0^{\rm th}$ order.

 The metric $g$ is a weak $L^{2,p}$ solution of the equation (2.5), 
i.e. (2.5) holds when it is paired with any $L_{o}^{2,p'}$ test form 
$h$, $(p^{-1} + (p' )^{-1} =$ 1) and integration by parts is performed 
twice; here $L_{o}^{2,p'}$ is the closure of the space of smooth 
functions of compact support in $M$ in the $L^{2,p'}$ norm.

  Given this setup, to control the boundary behavior it is necessary to 
make a specific choice of harmonic coordinates. Thus, given the boundary 
metric $\gamma \in L^{k+2, q}(\partial M)$, 
(or $C^{m+2, \alpha}(\partial M)$), choose local harmonic coordinates $u_{a}$, 
$a = 1,2,3$ for $\partial M$ w.r.t. $\gamma$. The coordinates $u_{a}$ are in 
$L^{k+3, q}(\partial M)$, (or $C^{m+3, \alpha}(\partial M)$), and 
$\gamma_{ab} = \gamma(\partial_{a}, \partial_{b}) \in L^{k+2, q}(\partial M)$, 
($C^{m+2,\alpha}(\partial M)$); here $\partial_{a} = 
\partial / \partial u_{a}$. Next, the coordinates $u_{a}$ may be extended to 
local harmonic coordinate functions for $M$ by solving a local Dirichlet 
problem: $\Delta_{g}u_{a} = 0$, with $u_{a}|_{\partial M}$ the given function 
$u_{a}$ on $\partial M$. Similarly, choose a local ``harmonic defining 
function'' $u_{4}$, i.e. $\Delta_{g}u_{4} = 0$, with 
$u_{4}|_{\partial M} = 0$. The metric 
$g = g_{ab} = g(\partial_{a}, \partial_{b})$ has optimal regularity in 
these coordinates.

  Clearly $g_{ab}|_{\partial M} = \gamma_{ab}$, for $a,b \leq 3$. However, 
$g_{4a}$ at $\partial M$ is not apriori determined by the boundary metric 
$\gamma$. The following Lemma shows that the components $g_{4a}$ satisfy
Neumann boundary conditions on $\partial M$.

\begin{lemma} \label{l 2.2}
Let $N = \nabla u_{4} / |\nabla u_{4}|$ be the unit normal at $\partial M$, 
$\nabla u_{4} = g^{4a}\partial_{a}$. Then the components $g^{4a}$ satisfy
the following Neumann boundary condition at $\partial M$:
\begin{equation} \label{e2.9}
N(g^{44}) = -6\lambda (g^{44})^{3/2} \  {\rm and} \ 
N(g^{4a}) = - \tfrac{1}{2}(g^{44})^{-1/2}g^{ab}\partial_{b}g^{44} 
-3\lambda \sqrt{g^{44}}g^{4a}, \ a < 4, \ {\rm on} \ \partial M,
\end{equation}
where $\lambda$ is given by $A = \lambda \cdot g$ at $\partial M$.
\end{lemma}

\noindent
{\bf Proof:}
Since $\Delta u_{a} = 0$ at $\partial M$, one has, at $\partial M$,
$$0 = <\nabla_{e_{i}}\nabla u_{a}, e_{i}> + 
<\nabla_{N}\nabla u_{\alpha}, N>,$$
where $e_{i}$ is an orthonormal basis for $\partial M$ at a given point. Write 
$\nabla u_{a} = (\nabla u_{a})^{T} + (\nabla u_{a})^{N}$. 
Then $<\nabla_{e_{i}}(\nabla u_{a})^{T}, e_{i}> = 0$, since $u_{a}$ 
is harmonic on $\partial M$, $a \leq 3$, and $u_{4} \equiv 0$ on 
$\partial M$. Further, $<\nabla_{e_{i}}(\nabla u_{a})^{N}, e_{i}> = 
A(e_{i}, e_{i})<\nabla u_{a}, N> = Hg^{4a} = 3\lambda g^{4a}$. This then gives 
\begin{equation} \label{e2.10}
<\nabla_{N}\nabla u_{a}, N> = -3\lambda g^{4a} \ \ {\rm at} \ \ \partial M.
\end{equation}

  Let $a = 4$. Since $\nabla u_{4} = g^{4a}\partial_{a}$, $|\nabla u_{4}|^{2} 
= g^{44}$, so $|\nabla u_{4}| = \sqrt{g^{44}}$. The term on the left in 
(2.10), with $a = 4$ is then
$$\frac{1}{|\nabla u_{4}|^{2}}<\nabla_{\nabla u_{4}}\nabla u_{4}, \nabla u_{4}> = \tfrac{1}{2}(g^{44})^{-1} <\nabla u_{4}, \nabla g^{44}> = \tfrac{1}{2}(g^{44})^{-1/2}N(g^{44}),$$
which, with (2.10), gives the first equation in (2.9).

  For the second equation, the left side of (2.10) may be written as
$$\frac{1}{|\nabla u_{4}|^{2}}<\nabla_{\nabla u_{4}}\nabla u_{a}, \nabla u_{4}> =  \frac{1}{|\nabla u_{4}|}N<\nabla u_{a}, \nabla u_{4}> - \frac{1}{|\nabla u_{4}|}<\nabla u_{a}, \nabla_{N}\nabla u_{4}>.$$
For the first term, compute
$$<\nabla u_{a}, \nabla u_{4}> = <g^{ac}\partial_{c}, g^{4b}\partial_{b}> 
= g^{ac}g^{4b}g_{bc} = g^{4a},$$ 
so that the first term is $\frac{1}{|\nabla u_{4}|}N(g^{4a})$.

  For the second term, using the gradient property, this is
$$\frac{1}{|\nabla u_{4}|}<N, \nabla_{\nabla u_{a}}\nabla u_{4}> = \frac{1}{2|\nabla u_{4}|^{2}}<\nabla u_{a}, \nabla g^{44}> = $$
$$\frac{1}{2|\nabla u_{4}|^{2}}<g^{ab}\partial_{b}, g^{cd} \frac{\partial g^{44}}{\partial u_{d}} \partial_{c}> = \frac{1}{2|\nabla u_{4}|^{2}}g^{ad}\partial_{d}g^{44}.$$
Combining these estimates gives the 2nd equation in (2.9).
{\endproof}

  We are now in position to deal with the proof of Proposition 2.1 itself. 
We begin with the lowest regularity situation, and so suppose $k = 1$ 
and $2 \leq q \leq p$. One has $R$ and $Ric$ in $L^{p}(M)$, and by 
assumption $s\in L^{1,q}(M)$ and $s|_{\partial M} \in L^{1,q}(\partial M)$. 
Hence, (2.5) and (2.7) give
\begin{equation} \label{e2.11}
\Delta (Ric_{ab}) \in  L^{-1,q}+L^{p/2} + L^{-1,p}, 
\end{equation}
on $M,$ where the right side of (2.11) denotes the sum of three terms, 
each in the respective spaces, cf. also [2, p.234]. (The $L^{-1,q}$ term 
comes from the $s$ terms on the right in (2.5), the $L^{p/2}$ term comes 
from the curvature terms in (2.5), while the $L^{-1,p}$ term comes from the 
$Q_{3}$ term in (2.7)). The coefficients of the Laplacian in (2.11) are in 
$L^{2,p}.$ We refer to [30] for treatment of Sobolev spaces with negative 
exponents, and recall that the dual space of $L^{-1,q}$ is $L_{o}^{1,q'}$. 

  A straightforward application of Sobolev embedding then shows that 
$L^{-1,p} \subset L^{-1,q}$ and $L^{p/2} \subset L^{-1,q}$ and so we 
may view the Laplacian in (2.11) as a mapping $\Delta: L^{1,q} 
\rightarrow  L^{- 1,q}$. Now elliptic boundary regularity theory for a 
Laplacian as in (2.8), cf. [30,Ch.2.7], [32,Thm.5.5.5'], shows that 
$Ric_{g} \in L^{1,q}$ provided the curvature $Ric_{g}$ is $L^{1,q}$ at the 
boundary $\partial M$.  By Lemma 1.3, the Ricci curvature of $g$ at 
$\partial M$ is determined algebraically by that of intrinsic metric 
$\gamma $, $s$ and $A$ at $\partial M$. By Lemma 1.1, (and passing to a new 
compactification equally as smooth as the original), one may assume that 
$A = 0$ at $\partial M$; alternately, it will be remarked below that the case 
$A \neq 0$ may be handled by the same arguments. It then follows from the 
assumptions on $Ric_{\gamma}$ and $s|_{\partial M}$ in Proposition 2.1 that 
$Ric_{g}|_{\partial M} \in L^{1,q}(\partial M)$ and hence 
$Ric_{g} \in L^{1,q}(M)$. 

  We now basically repeat this argument with (2.6) in place of (2.7). 
Thus, the left side of (2.6) is now in $L^{1,q}$, with Laplacian of the form 
(2.8) with $L^{2,p}$ coefficients. The lower order term $[Q_{1}]_{ab}$ is 
in $L^{1,p} \subset L^{1,q}$. For $i,j \leq 3$, since $g_{ij} = 
\gamma_{ij}$ at $\partial M$, it follows from elliptic boundary regularity 
that $g_{ij} \in L^{3,q}$, since the boundary metric 
$\gamma \in L^{3,q}(\partial M)$ by assumption. For the normal terms, 
$g^{4a}$, suppose as above, (w.l.o.g.), that $A = 0$ 
at $\partial M$ and work first with $g^{44}$. This satisfies the Neumann 
condition (2.9), where the coefficients of the vector field $N$ are in 
$L^{2,p}(M)$, and so in $C^{1,\alpha}(\partial M)$. Since $g^{44}$ also 
satisfies an equation of the form (2.6), with right-hand side of the form 
$Ric^{44} \in L^{1,q}(M)$, it follows from 
elliptic boundary regularity that $g^{44} \in L^{3,q}(M)$, c.f. 
[32,Thm.6.3.7] or [30,Ch.2.7.3] for instance. For the terms $g^{4a}$, 
$a < 4$, since now $\partial g^{44} \in L^{2,q}(M)$, the same arguments 
as above using the Neumann condition (2.9) on $g^{4a}$ give 
$g^{4a} \in L^{3,q}(M)$. 

  Thus we have $g_{ij} \in L^{3,q}(M)$, $i, j \leq 3$ and 
$g^{4a} \in L^{3,q}(M)$, $a \leq 4$. From this, it is an exercise in linear 
algebra to see that $g_{ab} \in L^{3,q}(M)$, $a, b \leq 4$; my thanks to 
Rafe Mazzeo for suggesting the argument below. Thus 
$g^{44} = (detg_{ab})^{-1}A_{44}$, where $A_{44}$ is the $(4,4)$ cofactor 
in the matrix $g_{ab}$. Since $A_{44}$ only involves $g_{ij}$, $i, j \leq 3$, 
the regularity above gives $A_{44} \in L^{3,q}(M)$, and hence 
$detg_{ab} \in L^{3,q}(M)$, since $g^{44} > 0$. The same reasoning on 
$g^{4a}$ then gives $A_{4a} \in L^{3,q}(M)$, $a \leq 3$. Now each 
determinant $A_{4a}$ may be expanded, (along its last row or column), to 
obtain a linear form in the variables $g_{4i}$, $i \leq 3$, with 
coefficients given by $2\times 2$ determinants. Thus one has a linear 
system of three equations in the three variables $g_{4i}$ with coefficients 
given by $2\times 2$ determinants. These $2\times 2$ determinants are 
cofactors in the $3\times 3$ matrix $g_{ij}$, $i, j \leq 3$. The determinant 
of the matrix associated to this $3\times 3$ linear system is then easily 
calculated to be just $-(detg_{ij})^{2}$. Since 
$detg_{ij}|_{\partial M} = det\gamma_{ij} > 0$, this linear system in 
invertible near $\partial M$. Hence, the variables $g_{4i}$ are rational 
expressions in $A_{4a}$, $g_{kl}$, $k,l \leq 3$ and $(detg_{kl})^{-1}$, each 
of which is in $L^{3,q}(M)$. It follows that $g \in L^{3,q}(M)$; this 
completes the proof in case $k = 1$ and $q \leq p$. 

  If $k = 1$ and $q > p$, then the work above gives $g \in L^{3,p}$ and 
by Sobolev embedding, $L^{3,p} \in L^{2,r}$, for any $r < \infty$, 
since $p > 4$. Choosing $r$ sufficiently large so that $q \leq r$, the 
work above with $r$ in place of $p$ gives $g \in L^{3,q}$, as required. 
This completes the proof in case $k = 1$.

  The proof for a given $k \geq 2$ follows exactly the same 2-step 
procedure, using induction from the regularity obtained at order $k-1$. 

  If $A \neq 0$ at $\partial M$, the same arguments are valid, with an extra 
bootstrap or iteration, since $A$ is of lower order.  Thus, for instance, 
working with (2.11), since $A^{2} \in L^{1,p}(M)$, $Ric_{g}|_{\partial M} 
\in L^{1,q} + L^{1-1/p,p}$ and so elliptic regularity gives 
$Ric_{g} \in L^{1-1/p,p}(M)$, (assuming $L^{1,q} \subset L^{1-1/p,p}$). 
In turn, this leads as before to $g \in L^{3-1/p,p}(M)$, which gives 
$A^{2} \in L^{2-1/p,p}(M)$. Using this new estimate for $A^{2}$ and iterating 
this process again leads to the same result as before.

 The proof of regularity in the case of H\"older spaces is 
also essentially the same although a little easier. Thus suppose $m = 1$, 
so that $\gamma \in C^{3,\alpha}$ and $s \in C^{1,\alpha} \subset L^{1,p}$, 
for all $p < \infty$. It follows that $D^{2}s \in L^{-1,p}$ and so the work on 
Sobolev space regularity above implies $g \in L^{3,p}$, for any $p < 
\infty$. The metric $g$ is thus an $L^{3,p}$ weak solution of (2.5) 
with $Ric|_{\partial M} \in C^{1,\alpha}$. Elliptic boundary regularity 
theory applied to (2.7), cf. [17,Ch.8], again implies that 
$Ric_{g} \in C^{1,\alpha}$ on $M$ and the same argument applied to (2.6) 
gives $g \in C^{3,\alpha}$. The proof when $m \geq 2$ follows in the same 
way, using the Schauder elliptic estimates. This then completes the proof 
in all cases.

 Regarding the estimate (2.4), the bound (2.2), together with bounds on 
the diameter and volume ratios of geodesic balls imply uniform 
$L^{2,p}$ bounds on the metric $\bar g$ in harmonic coordinates, as 
well as an upper bound on the number of such coordinate charts, cf. 
[2], [12]. Thus, the bound (2.4) follows from the fact that the 
elliptic regularity results above are effective, i.e. all the 
regularity statements are accompanied by estimates. 

 Finally, the equation (2.5) has real-analytic coefficients in the 
metric $\bar g,$ and smooth solutions of such equations with 
real-analytic boundary values are real-analytic, cf. [32, Ch.6.7]. 
This gives the statement on real-analyticity.

{\endproof}

\begin{remark} \label{r 2.3}
{\rm The method of proof above, and in particular Lemma 2.2, can also be 
used to prove other regularity results; for example regularity up to the 
boundary for metrics whose Ricci curvature is controlled up to the boundary.} 
\end{remark}

 Proposition 2.1 shows that the smoothness of the compactification 
$\bar g$ is determined by that of its scalar curvature $\bar s,$ and 
that of the boundary metric $\gamma $ on $\partial M.$ Since (2.1) is 
trace-free, one cannot improve this result eliminating the dependence 
on the scalar curvature $\bar s.$ An improvement can be obtained only 
by choice of a suitable representative of the conformal class, i.e. a 
suitable choice of gauge. From the formulas (1.4)-(1.5), a natural 
choice of gauge near $\partial M$ is that given by a geodesic defining 
function. However it seems to be difficult to prove higher order 
regularity directly in this gauge. The Yamabe, i.e. constant scalar 
curvature, gauge appears to be much better in this respect. 

 This leads to the following result, which is essentially Theorem 0.1.
\begin{theorem} \label{t 2.4.}
  Let $(M, g)$ be an AH Einstein 4-manifold which admits a $L^{2,p}$ 
conformal compactification $\bar g = \rho^{2}g$, $p > 4$. If, for a 
given $m \geq  3$ and $\alpha\in (0,1),$ the boundary metric $\gamma  = 
\bar g|_{\partial M}$ is $C^{m,\alpha}$ smooth, then there is another 
(possibly equal) conformal compactification $\widetilde g$ of g, with 
$\widetilde g|_{\partial M} = \gamma ,$ such that $\widetilde g$ is 
$C^{m,\alpha}$ smooth. If $\gamma $ is real-analytic, then there is a 
real-analytic compactification $\widetilde g.$

 Further, the estimate (2.4) holds for $\widetilde g$ without any 
dependence on the scalar curvature $\widetilde s.$
\end{theorem}

\noindent
{\bf Proof:}
Suppose $\gamma\in C^{m,\alpha}$, $m \geq  3$, and that $\bar g$ is an 
$L^{2,p}$ compactification of $g$. Let $\widetilde g$ be a constant 
scalar curvature metric conformal to $\bar g$ on $M$ with $\widetilde 
g|_{\partial M} = \gamma .$ Thus, for $\widetilde g = u^{2}\cdot \bar 
g,$ the function $u > 0$ is a solution to the Dirichlet problem for 
the Yamabe equation
\begin{equation} \label{e2.12}
u^{3}\mu  = - 6\bar \Delta u + \bar su, 
\end{equation}
on $M$, with $u \equiv  1$ on $\partial M$ and $\widetilde s = \mu$ = 
const. It is simplest to choose $\mu = -1$. Then standard methods in 
elliptic PDE give an $L^{2,p}$ solution to this Dirichlet problem, 
just as in the negative scalar curvature case of the Yamabe problem on compact 
manifolds. We refer to [29, Thm.1.1] and references therein, (cf. also 
[15,Remark]), for a proof, at least when $\bar g \in C^{2,\alpha}$. The 
same proof holds for $\bar g \in L^{2,p}$, $p > 4$; alternately, the 
compactness result of [22] implies that if $\bar g_{j} \in 
C^{2,\alpha}(M)$ converges to $\bar g$ in $L^{2,p}(M)$, then the Yamabe 
metrics $\widetilde g_{j}$ also converge in $L^{2,p}(M)$ to an 
$L^{2,p}$ Yamabe metric $\widetilde g \in [\bar g]$, with $u \in L^{2,p}(M)$. 

 Since the Bach equation (2.1) is conformally invariant, the metric 
$\widetilde g$ is an $L^{2,p}$ weak solution of the equation (2.5). 
Hence, since $\widetilde s$ is constant, Proposition 2.1 implies that 
$\widetilde g$ is as smooth as the boundary metric $\gamma$, which 
completes the proof.

 The estimate (2.4) follows as in the proof of Proposition 2.1, since 
the scalar curvature $\widetilde s$ is apriori controlled.
{\endproof}

 Theorem 2.4 gives the optimal regularity near $\partial M$ of a 
conformal compactification in terms of the regularity of the intrinsic 
metric $\gamma $ on $\partial M,$ assuming there is a $L^{2,p}$ 
conformal compactification of $g$, for some $p > $ 4. 

 We note also the following, essentially immediate, corollary.
\begin{corollary} \label{c 2.5.}
  Let $g$ be an AHE metric on M, which admits a $L^{2,p}$ conformal 
compactification, for some $p > $ 4, for which the boundary metric 
$\gamma $ is in $C^{m,\alpha}(\partial M),$ for some $m \geq  3$.

 Then the geodesic compactification $\bar g$ associated to $\gamma $ is 
at least a $C^{m-1,\alpha}$ compactification.
\end{corollary}

\noindent
{\bf Proof:}
 Let $\widetilde g$ be the $C^{m,\alpha}$ (Yamabe) compactification of 
$g$, given by Theorem 2.4. Then the geodesic compactification $\bar g$ 
and the Yamabe compactification $\widetilde g$ are related by a 
conformal factor $u$ satisfying (1.10). The coefficients and right side 
of the $1^{\rm st}$ order system (1.10) are in $C^{m,\alpha}$ and 
$C^{m-1,\alpha}$ respectively, and so the solution $u$ with $u \equiv $ 
1 on $\partial M$ is a $C^{m-1,\alpha}$ function on $\bar M.$ Hence 
$\bar g$ is also $C^{m-1,\alpha}$ on $\bar M.$
{\endproof}

\begin{remark} \label{r 2.6.}
  {\rm Theorem 2.4 does not hold in odd dimensions $n \geq $ 5. Namely 
by the result of Graham-Lee [19], there are $C^{\infty}$ metrics 
$\gamma $ on $S^{n-1}$ which are boundary metrics of $C^{n-2,\alpha}$ 
compactifications $\bar g$ of AHE metrics $g$ on the $n$-ball $B^{n}.$ 
However, for generic $\gamma ,$ such compactifications $\bar g$ have a 
non-zero $\rho^{n-1}\log \rho $ term in the asymptotic expansion (0.4) 
of $\bar g$ near $\partial M,$ cf. [18] for example. Hence, such 
metrics are at best only $C^{n-1,\alpha}$ smooth. It is unknown if 
Theorem 2.4 holds in even dimensions $n > 4$.}
\end{remark} 

 The geodesic compactification $\bar g$ is $C^{m-1,\alpha}$ only off 
the cutlocus $\bar C$ of $\partial M$ in $(\bar M, \bar g);$ at the 
cutlocus $\bar C,$ the metric $\bar g$ becomes singular, (although of 
course the Einstein metric $g$ is smooth). However, any smooth 
approximation to the geodesic defining function $t$ gives a smoothing 
to the compactification $\bar g.$ 

\smallskip

  We conclude this section with the following application of Corollary 2.5. 
First, let $\{g_{i}^{*}\}$ be a sequence of $L^{2,p}$ compactifications, 
$p > 4$, of AH Einstein metrics $(M, g_{i})$ with $C^{m,\alpha}$ boundary 
metrics $\gamma_{i}$, $m \geq 1$. Suppose (2.2) holds uniformly for 
$\{g_{i}^{*}\}$, 
and that the bounds on the diameter and volume ratios for geodesic balls 
of $g_{i}^{*}$ hold uniformly. Then a standard application of the $L^{p}$ 
Cheeger-Gromov compactness theorem, cf. [2], [12] and references therein, 
implies that $\{g_{i}^{*}\}$ is precompact, in that there is a subsequence 
converging in the weak $L^{2,p}$ and $C^{1,\alpha}$, 
$\alpha < 1-\frac{4}{p}$, topologies to an $L^{2,p}$ limit metric 
$g_{\infty}^{*}$ on $\bar M$.

  The following result shows that this convergence can be strengthened, 
given suitable control on the boundary metrics.
\begin{proposition} \label{p 2.7}
For $\{g_{i}^{*}\}$ as above, suppose the boundary metrics $\gamma_{i}$ are 
bounded in the $C^{m,\alpha}$ topology on $\partial M$, for some $m \geq 3$. 
Then, for the subsequence above, the $C^{m-1, \alpha}$ geodesic 
compactifications $\bar g_{i}$ determined by $\gamma_{i}$ converge, away 
from their cutlocus and in the $C^{m-1, \alpha'}$ topology, to the 
$C^{m-1, \alpha}$ limit $\bar g_{\infty}$, for any $\alpha' < \alpha$. 
Further, the distance of the cutlocus of each $\bar g_{i}$ to $\partial M$ 
is uniformly bounded below.
\end{proposition}

\noindent
{\bf Proof:} 
Corollary 2.5 and the associated (H\"older) bound (2.4) imply a uniform bound 
on $\{\bar g_{i}\}$ in the $C^{m-1, \alpha}$ topology on $M$, away from the 
cutlocus. Since $m \geq 3$, the curvature of $\bar g_{i}$ is thus uniformly 
bounded, which, by standard Riemannian geometry, implies a uniform lower 
bound on the distance of the cutlocus of $\bar g_{i}$ to $\partial M$. 
Given such a uniform bound on $\{\bar g_{i}\}$, it is then again standard 
that one has $C^{m-1, \alpha'}$ convergence to the $C^{m-1, \alpha}$ limit 
$\bar g_{\infty}$, for any $\alpha < \alpha'$; this is essentially the 
Arzela-Ascoli theorem in harmonic coordinates, c.f. [2], [12].
{\endproof}

  Using the $C^{1,\alpha}$ compactness, the geodesic compactifications 
$\bar g_{i}$ may be smoothed near the cutlocus to obtain $C^{m-1,\alpha'}$ 
convergence of the smoothed metrics on all of $\bar M$.

\section{Uniqueness.}
\setcounter{equation}{0}

 In this section, we prove the uniqueness theorem, Theorem 0.2. Let $g$ 
be an AH Einstein metric on a 4-manifold $M$, with $C^{m+1,\alpha}$ 
boundary metric $\gamma$, $m \geq  3$. By Corollary 2.5, we may assume 
that the geodesic compactification $\bar g$ associated with $\gamma $ 
is a $C^{m,\alpha}$ compactification, so that $\bar g$ has a 
Fefferman-Graham expansion
\begin{equation} \label{e3.1}
g_{t} = g_{(0)} + t^{2}g_{(2)} + t^{3}g_{(3)} + ... t^{m}g_{(m)}+ 
O(t^{m+\alpha}). 
\end{equation}
The coefficients are defined by
\begin{equation} \label{e3.2}
g_{(j)} = \tfrac{1}{j!}({\cal L}_{\bar \nabla t}^{(j)}\bar g)|_{\partial M}, 
\end{equation}
where ${\cal L}^{(j)}$ is the $j$-fold Lie derivative. Observe that 
although the expression (3.2) gives symmetric bilinear forms on 
$TM|_{\partial M},$ the vector $\bar \nabla t  \in $ Ker $g_{(j)}$ 
for all $j$ and so $g_{(j)}$ is uniquely determined by its restriction 
to $T(\partial M).$ Hence we view $g_{(j)}$ as bilinear forms on 
$\partial M.$

 The term $g_{(0)} = \gamma ,$ while the term $g_{(1)},$ equal to the 
$2^{\rm nd}$ fundamental form of $\partial M$ in $(M, \bar g),$ 
vanishes. Using the formulas (1.4)-(1.5) and (1.12)-(1.15), the term 
$g_{(2)}$ is given by
\begin{equation} \label{e3.3}
g_{(2)} = -\frac{1}{2}(Ric_{\gamma}-  \frac{s_{\gamma}}{4}\gamma ), 
\end{equation}
while the $g_{(3)}$ term satisfies
\begin{equation} \label{e3.4}
tr_{\gamma}g_{(3)} = 0, \delta_{\gamma}g_{(3)} = 0, 
\end{equation}
i.e. $g_{(3)}$ is transverse traceless. However, beyond the relations 
(3.3)-(3.4), the Einstein equations at $\partial M$ do not determine 
the coefficients $g_{(j)}, j \geq $ 3. In particular, the term 
$g_{(3)}$ is not apriori determined by the Einstein equations, for a 
given choice of boundary metric. These results follow from the work of 
Fefferman-Graham [16], cf. also [18], [25]. Related results hold in 
higher dimensions, given suitable boundary regularity, up to the order 
$g_{(n-1)}.$

\begin{remark} \label{r 3.1.}
  {\rm The term $g_{(3)}$ has the following interpretation from the 
AdS/CFT correspondence. First, the expansion (3.1) easily leads to an 
expansion for the volume of the geodesic 'spheres' $S(r) = \{x\in M: 
t(x) = 2e^{- r}\}$, of the form
\begin{equation} \label{e3.5}
vol S(r) = v_{(0)}e^{3r} + v_{(2)}e^{r} + O(e^{-\alpha r}), 
\end{equation}
cf. again [18] for instance. The coefficients $v_{(0)}$ and $v_{(2)}$ 
in (3.5) depend on the compactification $\bar g,$ and so are not 
invariantly attached to $(M, g)$. (The term $v_{(3)}$ vanishes by 
(3.4)). Let $B(r) = \{x\in M: t(x) \geq  2e^{- r}\}$ be the associated 
geodesic 'ball'. Then integrating (3.5) over $r$ gives
\begin{equation} \label{e3.6}
vol B(r) = \frac{1}{3}v_{(0)}e^{3r} + v_{(2)}e^{r} + V + o(1). 
\end{equation}
Now general reasoning from the AdS/CFT correspondence, cf. [38], leads 
to the conclusion that the constant term $V$ in (3.6) is in fact 
independent of the compactification $\bar g$ and depends only on $(M, 
g)$. 

 The term $V$ is the renormalized volume (or action, up to a 
multiplicative constant) of the AH Einstein metric $(M, g)$. In fact, 
$V$ may be computed invariantly in terms of the $L^{2}$ norm of the 
Weyl curvature $W$ of $(M, g)$ as
\begin{equation} \label{e3.7}
\frac{1}{8\pi^{2}}\int_{M}|W|^{2}dV = \chi (M) -  \frac{3}{4\pi^{2}}V. 
\end{equation}
cf. [3]. Note in particular that (3.7) implies $V \leq  \frac{4\pi^{2}}{3}\chi 
(M).$}

 {\rm Let $dV$ be the differential of $V$, acting on infinitesimal AH 
Einstein deformations $h$ of a given AH Einstein metric $(M, g)$, so 
that $g_{s} = g+sh$ is an AH Einstein metric to first order in $s$. Let 
$h_{(0)}$ be the induced first order variation of the boundary metric 
$\gamma_{s}$ at $\gamma .$ Then $dV$ is given by
\begin{equation} \label{e3.8}
dV_{g}(h) = -\frac{1}{4}\int_{\partial M}< g_{(3)}, h_{(0)}> dV_{\gamma}, 
\end{equation}
where the inner product and volume form are w.r.t. $\gamma$, c.f. again [3]. 
 Although (3.8) implies that $dV$ is determined by the behavior at 
$\partial M,$ $dV$ is not intrinsically determined by the boundary metric 
$\gamma$; it depends on the global AH Einstein filling $(M, g)$. 
A formula similar to (3.8) also holds in dimensions $n \geq $ 4, when 
$g_{(3)}$ is replaced by $g_{(n-1)}$ and lower order terms, cf. [3], 
and also [14], [33].}
\end{remark}

 For convenience, we restate Theorem 0.2 as follows. Define two 
manifolds $M^{1}$ and $M^{2}$ to be commensurable if $M^{1}$ and 
$M^{2}$ have covering spaces $\bar M^{1}, \bar M^{2}$ which are 
diffeomorphic. This is equivalent to the statement that the universal 
covers are diffeomorphic.

\begin{theorem} \label{t 3.2.}
  Let $(M, g)$ be an AH Einstein 4-manifold with $C^{7,\alpha}$ 
boundary metric. Then the data $(\gamma , g_{(3)})$ on $\partial M$ 
uniquely determine $(M, g)$ up to local isometry, i.e. if $g^{1}$ and 
$g^{2}$ are two such AH Einstein metrics on manifolds $M^{1}$ and 
$M^{2},$ with $\partial M^{1} = \partial M^{2} = \partial M$ such that, 
w.r.t. geodesic compactifications,
\begin{equation} \label{e3.9}
\gamma^{1} = \gamma^{2} \ \ {\rm and} \ \ g_{(3)}^{1} = g_{(3)}^{2}, 
\end{equation}
then $g^{1}$ and $g^{2}$ are locally isometric and the manifolds 
$M^{1}$ and $M^{2}$ are commensurable.
\end{theorem}

 The proof will be carried out in several steps below. The main issue 
is to prove that $g^{1}$ and $g^{2}$ are isometric within a collar 
neighborhood $U$ of $\partial M$ in $M$; given this it is 
straightforward to prove that $g^{1}$ and $g^{2}$ are then everywhere 
locally isometric. The basic idea to establish uniqueness within $U$ is 
to set up a suitable Cauchy problem for a conformal compactification 
within $U$, and then prove uniqueness of solutions to the Cauchy 
problem.

{\bf Step 1.}
  Let $g^{1}$ and $g^{2}$ be AH Einstein metrics on $M$ satisfying 
(3.9). By Corollary 2.5, the geodesic compactifications $\bar g^{1}$ 
and $\bar g^{2}$ of $g^{1}$ and $g^{2}$ are $C^{6,\alpha}$ 
compactifications. The discussion preceding Theorem 3.2 implies that 
the first 4 terms $g_{(j)},$ 0 $\leq  j \leq $ 3, of the Taylor 
expansion (3.1) for $g^{1}$ and $g^{2}$ agree, i.e.
\begin{equation} \label{e3.10}
g_{(j)}^{1} = g_{(j)}^{2}, j \leq  3. 
\end{equation}
However the geodesic defining functions $t^{i}$ for $g^{i}$ are not 
necessarily the same. We rectify this by means of a suitable 
diffeomorphism. Namely, a geodesic compactification (1.11) gives rise 
to a natural identification
$$I = I_{\bar g}: U \rightarrow  I\times \partial M, x \rightarrow  
(t(x), \sigma_{x}(0)), $$
where $\sigma_{x}$ is the unique $\bar g$ geodesic starting in 
$\partial M$ through $x$. For distinct AH Einstein metrics $g = g^{1}$ 
and $g^{2}$ with the same boundary metric $\gamma $ as above, the 
resulting identifications are distinct, although of course equivalent. 
Namely, in a possibly smaller collar neighborhood also called $U$, the 
diffeomorphism
$$\phi : U \rightarrow  U, \ \  \phi (I_{\bar g^{2}}^{-1}(t^{2}(x), 
\sigma_{x}^{2}(0))) = I_{\bar g^{1}}^{-1}(t^{1}(x), \sigma_{x}^{1}(0)), 
$$
has the effect that $\phi^{*}t^{2} = t^{1}: U \rightarrow  {\Bbb R}$. 
Observe that $\phi \in C^{6, \alpha}$; this is because the vector fields 
$\nabla_{g^{1}}t^{1}$ and $\nabla_{g^{2}}t^{2}$ are $C^{6, \alpha}$ and 
so have $C^{6, \alpha}$ flows. The map $\phi$ is a composition of these 
two flow maps.

  Thus, given the fixed metric $\bar g = \bar g^{1},$ we pull back the 
metric $\bar g^{2}$ to the metric
$$\bar {\bar g}^{2} \equiv  \phi^{*}\bar g^{2}. $$
The $C^{5, \alpha}$ metric $\bar {\bar g}^{2}$ is of course isometric to 
$\bar g^{2}$ in $U$, has geodesic defining function $t^{1},$ and hence a 
splitting w.r.t. $t^{1}.$ Further, the metrics $\bar {\bar g}^{2}$ and $\bar 
g^{1}$ have the same boundary metric $\gamma$.

  Since $\phi = id$ on $\partial M$, the terms $\bar g_{j}^{2}$ and 
$\bar {\bar g}_{j}^{2}$ are equal for $j \leq 3$; this can be seen 
directly from the expressions (3.3) and (3.8), cf. also [14] and referemces 
therein. It follows then that the Taylor expansions of $\bar g^{1}$ and 
$\bar {\bar g}^{2}$ w.r.t. $t^{1}$ agree up to order 3. In the following, 
we will always assume that $\bar g^{2}$ is pulled back in this way to make 
it comparable to a given $\bar g.$

{\bf Step 2.}
 As discussed in \S 2, Einstein metrics $g$ and their compactifications 
$\bar g$ are solutions of the conformally invariant Bach equation (2.5) 
in dimension 4,
\begin{equation} \label{e3.11}
2D^{*}DRic +\frac{2}{3}D^{2}s +\frac{1}{3}\Delta s \cdot g +{\cal 
R}_{1} = 0. 
\end{equation}
Here and below, we will usually drop the overbar from the notation.

 Because of the conformal as well as diffeomorphism invariance of 
(3.11), one must choose suitable gauges, i.e. representatives of the 
conformal and diffeomorphism actions, in order to prove any uniqueness. 
For the conformal gauge, we choose the geodesic compactification, while 
for the diffeomorphism gauge, we use harmonic coordinates.

 The AH Einstein metric $(M, g)$ has a $C^{5,\alpha}$ geodesic 
compactification $\bar g$ with boundary metric $\gamma ,$ and with 
$C^{6,\alpha}$ geodesic defining function $t$. By (1.5), the scalar 
curvature $s = \bar s$ is given by $s = - 6\frac{\Delta t}{t},$ and in 
local harmonic charts for a neighborhood $U$ of $\partial M$ one thus 
has
$$s = - 6g^{ij}t^{-1}\partial_{i}\partial_{j}t.$$
It follows that the terms $D^{2}s$ and $\Delta s$ in (3.11) involve at 
most the second derivatives of the metric tensor $\bar g_{ij}$, with 
coefficients that are at least $C^{1,\alpha}.$ 
Thus, (3.11) may be rewritten as
\begin{equation} \label{e3.12}
D^{*}DRic + Q_{2}(x, g, \partial_{j}g) = 0, 
\end{equation}
where $Q_{2}$ involves $g$ and its derivatives only up to order 2, with all 
coefficients at least $C^{1,\alpha}$. As in (2.6)-(2.7), one may then 
rewrite the system (3.11) in a harmonic coordinate atlas ${\cal A} $ 
covering $\partial M$ as
\begin{equation} \label{e3.13}
\Delta\Delta g + Q_{3}(x, g, \partial_{j}g) = 0, 
\end{equation}
where $Q_{3}$ involves derivatives of $g$ only up to order 3.

 By Step 1, if $g^{1}$ and $g^{2}$ are two distinct metrics satisfying 
(3.9), we may assume that the associated geodesic defining functions 
are the same.

{\bf Step 3.}
 In this step, we set up and prove uniqueness for the Cauchy problem 
for the system (3.13).

 From the work in Steps 1 and 2, to each AH Einstein metric $g$ with 
boundary metric $\gamma ,$ we have associated a geodesic 
compactifiction $\bar g$ defined in a collar neighborhood $U$ of 
$\partial M,$ with fixed defining function $t$. In a harmonic atlas for 
$\bar U,$ the metric $\bar g$ satisfies the system (3.13). This is a 
$4^{\rm th}$ order, non-linear, elliptic system in the metric $\bar g.$ 
Further, in these local coordinates, the Cauchy data on $\partial M$ 
takes the form
\begin{equation} \label{e3.14}
\partial_{t}^{(j)}(g_{ij}) = g_{(j)}, 0 \leq  3 \leq  j, \ \ {\rm on} \ 
\  \partial M. 
\end{equation}
Clearly, $\partial M$ is non-characteristic for the Cauchy problem 
(3.13)-(3.14). As explained at the beginning of \S 3, the data (3.14) 
are determined by the data (3.9).

 We now claim that this coordinate Cauchy problem has a unique solution 
in a possibly smaller neighborhood $U'  \subset  U$ of $\partial M.$ 
Given this for the moment, if $g^{2}$ is another AH Einstein metric 
with boundary metric $\gamma ,$ then by construction in Steps 1 and 2, 
the geodesic compactification $\bar g^{2}$ is also a solution, (in 
local harmonic coordinate charts), to the Cauchy problem (3.13) with 
the same boundary Cauchy data (3.14). Hence, uniqueness to the 
coordinate Cauchy problem implies that the metrics $g^{1}$ and $g^{2}$ 
are isometric in $U$.

 With regard to uniqueness of the coordinate Cauchy problem, first note 
that the symbol (or characteristic polynomial) of $\Delta $ is
$$\sigma (\Delta ) = |\xi|^{2} = g^{ij}\xi_{i}\xi_{j}: T^{*}(M) 
\rightarrow  {\Bbb R} , $$
where $g^{ij}$ is the metric induced on the cotangent bundle. Here and 
below, all computations are w.r.t. the compactification $\bar g,$ but 
we omit the overbar from the notation. Hence, the leading order symbol 
of the Bach equation in the form (3.13) is
\begin{equation} \label{e3.15}
\sigma (B) = \sigma (\Delta\Delta ) = |\xi|^{4}. 
\end{equation}
This symbol has of course no real characteristics. However, it does 
have double complex characteristics. Namely, for $\xi\in T^{*}(\partial 
M)$ with $|\xi| =$ 1, the roots of the characteristic form of the 
leading term $\Delta\Delta$,
$$p(x, \xi , \tau ) = \sigma (\xi +\tau dt) = 0, $$
are given by
$$\tau  = \pm i, $$
independent of $x\in U$ and $\xi $ in the unit sphere bundle within 
$T^{*}(\partial M).$ Thus, the operator $\Delta\Delta $ has constant, 
pure imaginary, double characteristics.

 The proof of uniqueness for this coordinate Cauchy problem now 
essentially follows from the Calder\'on uniqueness theorem, cf. [10] 
and especially [11, Theorem 11]. A clear exposition of this result is 
also given by Nirenberg in [34, \S 6-\S 7], (however only in the case 
of linear equations with $C^{\infty}$ coefficients); cf. also [36]. We 
describe below how to reduce the uniqueness problem to the class of 
problems solved in [11, Theorem 11].

 The main issue is to reduce the non-linear Cauchy problem to a linear 
one. Note that while the nonlinearity in the lower order term $Q_{3}$ 
in (3.13) is complicated, the nonlinearity in the leading order term 
just comes from the fact that the "unknown" $g$ enters in the Laplacian 
$\Delta ,$ of the form (2.8). 

 We do this following the elegant method of [10, \S 5]. Thus, the system 
(3.13) is a nonlinear system of 10 $4^{\rm th}$ order PDE's in 10 
unknowns $g = g_{ij} = g_{ji}$ on (a domain in) $({\Bbb R}^{4})^{+}.$ 
Let $u$ denote a variable vector in ${\Bbb R}^{10},$ (so that $u$ 
corresponds to the metric $g$), and $\{u_{\alpha}\}$ the collection of 
all partial derivatives of $u$ of order $\leq $ 4. The system (3.13) 
may then be formally expressed as
\begin{equation} \label{e3.16}
F(x, u, u_{\alpha}) = 0, 
\end{equation}
where $F: ({\Bbb R}^{4})^{+}\times {\Bbb R}^{10}\times {\Bbb R}^{64} 
\rightarrow  {\Bbb R}^{10}$ is $C^{1,\alpha}$ smooth. Write
\begin{equation} \label{e3.17}
F(x, u, u_{\alpha}) = \Delta_{u}\Delta_{u}u + F_{3}(x, u, u_{\alpha}), 
\end{equation}
where $F_{3}$ corresponds to the term $Q_{3}$ in (3.13), so that 
$F_{3}$ has order 3.

 Now suppose $u$ and $v$ are two solutions of (3.16), corresponding to 
metrics $g = g^{1}$ and $g^{2}$, with say $v$ fixed. One then has
$$0 = \Delta_{u}\Delta_{u}u -  \Delta_{v}\Delta_{v}v + F_{3}(x, u, 
u_{\alpha}) -  F_{3}(x, v, v_{\alpha}) = \Delta_{u}\Delta_{u}(u- v) + 
H(x, u, u_{\alpha}) -  H(x, v, v_{\alpha}), $$
where
$$H(x, u, u_{\alpha}) = \Delta_{u}\Delta_{u}v +  F_{3}(x, u, 
u_{\alpha}). $$
Note that in terms of the metrics $u = g = g_{1}, v = g_{2},$ 
(subscripted here for convenience),
$$\Delta_{u}\Delta_{u}v = g_{1}^{ij}g_{1}^{kl}\partial_{ijkl}g_{2}. $$

 Now as in [10, \S 5] the mean value theorem applied to $H$, (with $x$, 
$v$ fixed and $u$ varying), gives
\begin{equation} \label{e3.18}
H(x, u, u_{\alpha}) -  H(x, v, v_{\alpha}) = (u-v)\int_{0}^{1}H_{u}[x, 
v+(u- v)s, v_{\alpha}+(u_{\alpha}- v_{\alpha})s]ds +
\end{equation}
$$+ \sum_{\alpha}(u_{\alpha}- v_{\alpha})\int_{0}^{1}H_{\alpha}[x, 
v+(u- v)s, v_{\alpha}+(u_{\alpha}- v_{\alpha})s]ds.$$
Substitute the solutions $u = u(x)$ and $v = v(x)$ in all terms inside 
the integrals in (3.18), so that the  integrals then become  
coefficient functions in $x$. Hence, (3.18) becomes a $3^{\rm rd}$ 
order linear system in the unknown $u- v.$ It follows that one has a 
solution $u- v$ of the linear $4^{\rm th}$ order system
\begin{equation} \label{e3.19}
\Delta_{u}\Delta_{u}(u- v) + H(x, (u- v), (u_{\alpha}- v_{\alpha})) = 
0. 
\end{equation}

 The leading order symbol of (3.19) is given by (3.15), and $u- v$ has 
0 Cauchy data on $\partial M.$ Hence, if the leading coefficients in 
(3.19) are in $C^{1,\beta}, \beta  > $ 0, and the lower order 
coefficients are bounded and measurable, then [11, Theorem 11] implies 
that $u = v$ in a neighborhood $U$ of $\partial M.$ We have assumed 
that $u = g^{1}, v = g^{2}\in C^{5,\alpha}.$ This implies that the 
lower order coefficients are at least in $C^{\alpha},$ while the 
leading order coefficient is in $C^{5,\alpha}.$

 This completes the proof of uniqueness within a collar neighborhood 
$U$. The last step is to extend this to the filling manifolds $M^{1}$ 
and $M^{2}$ of $g^{1}$ and $g^{2}.$

{\bf Step 4.}
 Suppose $g^{1}$ and $g^{2}$ are two AHE metrics on manifolds $M^{1}$ 
and $M^{2}$ which agree, up to diffeomorphism on a collar neighborhood 
$U$ of $\partial M = \partial M^{i},$ i.e. there is a diffeomorphism 
$\phi : U \rightarrow U, \phi  = id$ on $\partial M,$ such that
\begin{equation} \label{e3.20}
\phi^{*}g^{2} = g^{1}. 
\end{equation}

 We claim that $g^{1}$ and $g^{2}$ are locally isometric, i.e. for all 
$x^{1}\in M^{1}$ there exists $x^{2}\in M^{2}$ together with small open 
balls $V^{i}\in M^{i}, x^{i}\in V^{i}$ and a diffeomorphism $\psi : 
V^{1} \rightarrow  V^{2},$ such that $\psi^{*}g^{2} = g^{1}$ on 
$V^{1}.$ To see this, let $K^{i}$ be a domain with compact closure in 
$M^{i}$ such that $\partial K^{i} \subset U$, so that $M^{i} = K^{i} 
\cup U$. For each $g^{i},$ we may cover $K^{i}$ by a finite collection 
of charts which are harmonic w.r.t. $g^{i},$ i.e. let ${\cal A}^{i}$ be 
a finite harmonic atlas for (a thickening of) $K^{i}$ w.r.t. $g^{i}.$ 
By (3.20), without loss of generality we may assume that the charts in 
${\cal A}^{1}$ restricted to $U \cap K^{1}$ are $\phi$-pullbacks of 
charts in ${\cal A}^{2}$ restricted to $U \cap K^{2}$. 

 Now it is well-known that in harmonic charts an Einstein metric is 
real-analytic and hence satisfies unique continuation. Thus, given the 
expression for the local components of $g^{1}$ in one local harmonic 
chart of ${\cal A}^{1},$ the expression for $g^{1}$ in all of the other 
finitely many harmonic charts of ${\cal A}^{1}$ is uniquely determined, 
by analytic continuation along paths. The same holds w.r.t. $g^{2}.$ 

 Thus, given $x^{1}\in K^{1},$ let $\sigma^{1}$ be an analytic path in 
$K^{1}$ joining $x^{1}$ to a point $x^{o} \in U \cap K^{1}$. Using the 
identification $\phi : U \rightarrow  U$ near $\partial M$ and 
analyticity, $\sigma^{1}$ gives rise to a unique path $\sigma^{2}$ in 
$K^{2},$ ending at a point $x^{2}\in K^{2}.$ Since $g^{2}$ and $g^{1}$ 
are isometric in $U$, analytic continuation along $\sigma^{1}$ and 
$\sigma^{2}$ implies that $g^{2}$ are $g^{1}$ are locally isometric 
near $x^{2}$ and $x^{1}.$ Of course, the local isometry may depend on 
the homotopy class of the path $\sigma^{1}.$ An alternate, but 
essentially similar argument is to show that the set of points where 
$g^{2}$ and $g^{1}$ are locally isometric is both open and closed, cf. 
also [27, Ch.6.6].

 Finally, since $g^{1}$ and $g^{2}$ are locally isometric, it follows 
that they are isometric in the universal covers of each $M^{i},$ and 
hence the manifolds $M^{1}$ and $M^{2}$ are commensurable.

{\endproof}

\begin{remark} \label{r 3.3.} {\bf (i).}
  {\rm We point out that the uniqueness, within a collar neighborhood, 
of self-dual AH Einstein metrics with real-analytic compactifications, 
has been proved by LeBrun in [28], using twistor methods.}

{\bf (ii).}
 {\rm The proof of Theorem 3.2 strongly uses the fact that dim $M =$ 4, 
since the conformally invariant Bach equation can be used in that 
situation. It is unknown if an analogous result holds in $n$-dimensions, 
i.e. whether the coefficients $(g_{(0)}, g_{(n-1)})$ uniquely determine 
an AH Einstein metric up to local isometry. Without working in the 
compactified setting, this would require a uniqueness result for the 
Cauchy problem for a highly degenerate elliptic system.}

{\bf (iii).}
 {\rm It follows of course from Theorem 3.2 that all the higher order 
terms $g_{(j)}$ in the Fefferman-Graham expansion (3.1) are uniquely 
determined by the pair $(\gamma , g_{(3)}).$ This also follows directly 
from an obvious analysis of the Bach equation at the boundary $\partial 
M.$}

{\bf (iv).}
 {\rm The hypothesis $\gamma\in C^{7,\alpha}$ is needed only for 
technical reasons arising from the proof. In the sequel paper [4], 
methods will be developed allowing one to use approximation arguments, 
so that the hypothesis $\gamma\in C^{7,\alpha}$ can be relaxed to 
$\gamma\in C^{3,\alpha}.$}
\end{remark}

 Theorem 3.2 implies that the isometry type of $(M, g)$ is determined 
by $(\gamma , g_{(3)})$ and the action of $\pi_{1}(M)$ on the universal 
cover, i.e. the representation of $\pi_{1}(M)$ as a subgroup of the 
isometry group Isom$(\widetilde M)$ of $\widetilde M.$ The examples 
constructed in \S 4.4 below are locally isometric, non-isometric 
metrics on a fixed manifold, with a fixed $(\gamma , g_{(3)}),$ but 
varying representation of $\pi_{1}(M).$

\section{Non-Uniqueness.}
\setcounter{equation}{0}

 In this section, we examine in detail several classes of examples 
which show that in general an AH Einstein metric is not uniquely 
determined by its conformal infinity. These examples will also 
illustrate the sharpness of the uniqueness result, Theorem 3.2. These 
classes of examples are AdS black hole metrics and are discussed in 
some detail in the literature on the AdS/CFT correspondence. cf. [23], 
[24], [33], [38], and also [6, (9.118)]. The black hole topologies may 
be arbitrary surfaces, i.e. $S^{2}, T^{2}$ or $\Sigma_{g},$ where 
$\Sigma_{g}$ is any oriented surface of genus $g \geq $ 2. The most 
interesting cases, (for the present purposes), are those of $S^{2}$ and 
$T^{2},$ which we treat first and last.

\medskip

{\bf \S 4.1.}
 We begin with a discussion of the AdS-Schwarzschild metric, following 
[23]. On the manifold $M = {\Bbb R}^{2}\times S^{2},$ consider the 
metric
\begin{equation} \label{e4.1}
g_{m} = g_{m}^{(+1)}= V^{-1}dr^{2} + Vd\theta^{2} + r^{2}g_{S^{2}(1)}, 
\end{equation}
where $g_{S^{2}(1)}$ is the standard metric of curvature $+1$ on 
$S^{2}$ and $V = V_{m}(r)$ is the function
\begin{equation} \label{e4.2}
V = 1+r^{2} -  \frac{2m}{r}. 
\end{equation}
The mass $m$ is any positive number, $m >$ 0; if $m < $ 0, the metric 
(4.1) has a singularity at $r =$ 0 and so it is no longer complete. The 
parameter $r$ runs over the interval $[r_{+}, \infty ),$ where $r_{+}$ 
is the largest root of the equation $V(r) = 0$. The locus $\Sigma  = 
\{r = r_{+}\}$ in $M$ is thus a totally geodesic round 2-sphere, of 
radius $r_{+}.$ The circular parameter $\theta$ runs over an interval 
$[0, \beta]$ of length $\beta$. Smoothness of the metric $g_{m}$ at 
$\Sigma$ requires that 
$$\lim_{r \rightarrow r^{+}} V^{1/2}\frac{d(V^{1/2})}{dr}\beta  = 2\pi;$$
otherwise, the metric has a cone singularity along and normal to $\Sigma$. 
It follows easily from this and (4.2) that $g_{m}$ is smooth everywhere 
exactly when
\begin{equation} \label{e4.3}
\beta  = \frac{4\pi r_{+}}{1+3r_{+}^{2}}. 
\end{equation}
Observe also that the radius $r_{+}$ increases monotonically from 
0 to $\infty $ as the mass parameter $m$ increases from 0 to $\infty .$

  If one sets $m = 0$ in (4.2) and $\beta = \infty$, then the metric (4.1) 
is the hyperbolic metric $H^{4}(-1)$ on the 4-ball $B^{4}$, (decomposed 
along equidistants from $H^{3}(-1) \subset H^{4}(-1)$). This can be seen 
by the change of coordinates $r = \sinh s$. Here of course the sphere $\Sigma$ 
has collapsed to a point. However, the metrics $g_{m}$ do not converge 
(globally) to the hyperbolic metric as $m \rightarrow 0$, due to the 
restriction on $\beta$ in (4.3). As $m \rightarrow 0$, $r_{+} \rightarrow 0$, 
and so $\beta \rightarrow 0$. Nevertheless, for $r$ large, the term 
$2m/r$ in (4.2) is small and so the local geometry of the metric $g_{m}$, 
for any $m > 0$, approximates hyperbolic geometry. In fact, it is easily 
verified that the metrics $g_{m}$ are conformally compact, with conformal 
infinity given by the conformal class of the product metric 
$S^{1}(\beta) \times S^{2}(1)$.
  
 The 1-parameter family of metrics $g_{m}$ are Einstein metrics 
satisfying (0.1), and are isometrically distinct, i.e. $g_{m_{1}}$ is 
not isometric to $g_{m_{2}}$ for $m_{1} \neq  m_{2}.$ The parameter 
$\beta $ in (4.3) however does not increase monotonically with $m$ 
or $r_{+}.$ In fact, $\beta $ has a maximal value $\beta_{o},$
$$\beta_{o} = 2\pi /\sqrt {3}, $$
achieved at $r_{+} = 1/\sqrt {3}$. As $m \rightarrow 0$, or 
$m \rightarrow  \infty$, $\beta  \rightarrow 0$. In particular, for any 
$m_{1} \neq  2/(3)^{3/2},$ there is an $m_{2} \neq  m_{1}$ such that 
the AH Einstein metrics $g_{m_{1}}$ and $g_{m_{2}}$ on $S^{2}\times 
{\Bbb R}^{2}$ are not isometric but have the same conformal infinity. 
This is the first example of non-uniqueness.

 As indicated above, as $\beta  \rightarrow 0$, these metrics degenerate, 
as does the conformal structure of the boundary metric. Observe also that 
since $\beta  \leq  \beta_{o},$ the boundary metrics $S^{1}(\beta )\times 
S^{2}(1)$ for $\beta  >  \beta_{o}$ are not achieved in this family.

\begin{remark} \label{r 4.1} {\bf (i).}
  {\rm There is another AH Einstein metric with conformal infinity 
$S^{1}(\beta )\times S^{2}(1).$ Namely let $\gamma $ be a geodesic in 
the hyperbolic space $H^{4}(- 1)$ and let $(M, g)$ $= (H^{4}(- 1)/{\Bbb 
Z} , g_{-1}),$ where the ${\Bbb Z} $ action is generated by translation 
of length $\beta $ along $\gamma .$ This hyperbolic metric also has 
conformal infinity given by $S^{1}(\beta )\times S^{2}(1).$ Note that 
the topological type here, ${\Bbb R}^{3}\times S^{1},$ is distinct from 
that of the Schwarzschild family. 

 In this situation, all values of the length $\beta $ may be realized 
as boundary metrics. Further, if one replaces the (pure hyperbolic) 
translation along $\gamma $ by a loxodromic translation, i.e. 
translation along $\gamma $ together with a rotation in the orthogonal 
$H^{3}(- 1),$ then the resulting conformal structure at infinity is a 
bent product $S^{1}(\beta )\times _{\alpha}S^{2}(1),$ where the angle 
$\alpha $ between the factors corresponds to the twist rotation.}

{\bf (ii).}
 {\rm There are a number of other explicit examples of $S^{2}$ black 
hole AdS metrics; for example the AdS Taub-Bolt metrics on non-trivial 
line bundles over $S^{2},$ cf. [24] and references therein.}
\end{remark}

{\bf \S 4.2.}
 Next, consider the class of AdS black hole metrics on surfaces $\Sigma 
 = \Sigma_{g},$ of genus $\geq $ 2. As above, on the manifold $M = 
{\Bbb R}^{2}\times \Sigma ,$ consider the metric
\begin{equation} \label{e4.4}
g_{m} = g_{m}^{(- 1)}= V^{-1}dr^{2} + Vd\theta^{2} + r^{2}g_{\Sigma}, 
\end{equation}
where $g_{\Sigma}$ is a hyperbolic metric on $\Sigma ,$ i.e. any point 
in the moduli space of Riemann surfaces. Now $V = V(r)$ is given by
\begin{equation} \label{e4.5}
V = - 1+r^{2} -  \frac{2m}{r}, 
\end{equation}
with $r \geq  r_{+},$ the largest root of $V(r) = 0$. As before, the 
locus $\Sigma  = \{r = r_{+}\}$ in $M$ is totally geodesic and 
isometric to $(\Sigma , g_{\Sigma}),$ and smoothness at the horizon 
requires $\theta\in [0,\beta )$ with
\begin{equation} \label{e4.6}
\beta  = \frac{4\pi r_{+}}{- 1+3r_{+}^{2}}. 
\end{equation}
The metrics $g_{m}$ are AH Einstein metrics on $M$, with conformal 
infinity $S^{1}(\beta )\times \Sigma ,$ and are non-isometric for 
distinct values of $m$.

 In contrast to the case of $S^{2},$ the function $\beta $ here is 
monotone decreasing as $m$ or $r_{+}$ increases, so that $\beta $ is a 
single valued function of $m$ or $r_{+}.$ Further, the metric $g_{m}$ 
is well-defined whenever $r_{+} >  1/\sqrt {3}$, which is equivalent to
$$m >  m_{o} = -  3^{-3/2}. $$
Hence the mass parameter may assume (some) negative values.

 When $m =$ 0, so that $r_{+} =$ 1, the metric $g_{o}$ is the 
hyperbolic metric on ${\Bbb R}^{2}\times \Sigma ,$ and $\beta $ has the 
value $2\pi .$ When $m \rightarrow  \infty , \beta  \rightarrow $ 0, 
while when $m \rightarrow  m_{o}, \beta  \rightarrow  \infty .$ Thus, 
at these extremes, both the metrics and the conformal infinity 
degenerate. In particular, we see that this family does not provide 
examples of non-uniqueness.

\medskip

{\bf \S 4.3.}
 Before proceeding to discuss $T^{2}$ AdS black hole metrics, in this 
subsection we review the well-known theory of Dehn surgery on 
hyperbolic 3-manifolds. This review mainly motivates the construction 
to follow in dimension 4 in \S 4.4, but also shows that uniqueness 
fails even in the category of conformally compact hyperbolic 
3-manifolds.

 Let $(T^{2}, g_{o})$ be a torus with a fixed flat metric $g_{o},$ 
representing a fixed point in the moduli space of flat structures on 
$T^{2}.$ Let $\sigma $ be a given simple closed geodesic in $T^{2},$ 
with length $L = L(\sigma ).$

 Next, let $\gamma $ be a complete geodesic in $H^{3}(-1)$ and let 
$T(R)$ be the $R$-tubular neighborhood about $\gamma $ in $H^{3}(-1).$ 
The metric on $T(R)$ then has the form
$$g_{-1} = dr^{2}+\sinh^{2}rd\theta^{2} + \cosh^{2}rds^{2}, $$
where $s$ is the parameter for $\gamma $ and $\theta\in [0,2\pi ].$ The 
boundary $\partial T(R)$ is a flat cylinder $S^{1}\times {\Bbb R} ,$ 
with metric
\begin{equation} \label{e4.7}
\widetilde g_{o} = \sinh^{2}Rd\theta^{2} + \cosh^{2}Rds^{2}. 
\end{equation}
Now choose $R$ so that 
$$2\pi \sinh R = L(\sigma ). $$
There is then a unique free ${\Bbb Z}$-action on the cylinder $\partial 
T(R)$ such that the quotient $S^{1}\times _{{\Bbb Z}}{\Bbb R} $ with 
the induced metric is the given flat torus $(T^{2}, g_{o})$ and such 
that the meridian circle $S^{1} = \partial D^{2}$ of length 
$2\pi\sinh R$ in the cylinder is mapped to $\sigma .$ 

 This action extends to an isometric action on $T(R)$ and so produces a 
hyperbolic metric $g_{- 1}$ on the solid torus $D^{2}\times S^{1},$ 
with boundary isometric to $(T^{2}, g_{o}),$ and with the geodesic 
$\sigma $ in $T^{2}$ bounding the disc $D^{2}$ in $D^{2}\times S^{1}.$ 
This metric is the tube of radius $R$ about the core closed geodesic 
$\gamma .$ Observe that the length of the core geodesic $\gamma ,$ of 
distance $R$ to $\partial T(R),$ is on the order of $O(\sinh^{-1}R) 
<< 1$, for $R$ large.

 It is clear that this hyperbolic metric extends to a complete 
hyperbolic metric on $D^{2}\times S^{1}$ with smooth conformal 
infinity. Since $(T^{2}, g_{o})$ is the metric $\widetilde g_{o}$ on 
$S^{1}\times {\Bbb R} $ divided out by the ${\Bbb Z} $ action, the 
conformal infinity is given by the conformal class $(T^{2}, 
[g_{\infty}])$ where
\begin{equation} \label{e4.8}
g_{\infty} = (e^{2R}d\theta^{2} + e^{2R}ds^{2})/{\Bbb Z} . 
\end{equation}

 The classes $[g_{\infty}]$ and $[g_{o}]$ do not agree, (although 
$[g_{\infty}] \rightarrow  [g_{o}]$ on any sequence where $L(\sigma ) 
\rightarrow\infty ).$ However, the construction above can easily be 
modified so that the conformal infinity is fixed instead of fixing the 
conformal structure $g_{o}$ on $\partial T(R).$ Namely, for any fixed 
$R$, write $s'  = s' (R) = \frac{\cosh R}{\sinh R}s,$ so that in these new 
coordinates, the metric $\widetilde g_{o}$ in (4.7) has the form
\begin{equation} \label{e4.9}
\widetilde g_{o} = \sinh^{2} R(d\theta^{2} + (ds' )^{2}). 
\end{equation}
Now divide $D^{2}(R)\times {\Bbb R} $ and $\partial D^{2}(R)\times 
{\Bbb R} $ by the same ${\Bbb Z} $ action as before, but with respect 
to the parameters $(\theta , s' )$ in place of $(\theta , s)$. This 
gives a complete hyperbolic metric on $D^{2}\times S^{1}$ with 
prescribed conformal infinity $(T^{2}, g_{o}),$ for any choice of 
closed geodesic $\sigma  \subset  (T^{2}, g_{o}).$

 Summarizing, the discussion above proves:

\begin{proposition} \label{p 4.2.}
  For any given flat structure $g_{o}$ on the torus $T^{2},$ and for 
any given simple closed geodesic $\sigma $ in $(T^{2}, g_{o}),$ there 
is a unique complete hyperbolic metric $g_{- 1}$ on the solid torus 
$D^{2}\times S^{1},$ with $(T^{2}, g_{o})$ as conformal infinity.
\end{proposition}
{\endproof}

 As $\sigma $ varies over the class of simple closed geodesics on 
$T^{2},$ the resulting hyperbolic metrics, although of course locally 
isometric, are not isometric since for instance the lengths of the core 
geodesics are distinct; compare with the discussion at the end of \S 3. 
In particular, there are infinitely many distinct hyperbolic 
3-manifolds, all diffeomorphic to $D^{2}\times S^{1},$ whose conformal 
infinity is an arbitrary but fixed $(T^{2}, g_{o}).$

 As $L = L(\sigma ) \rightarrow  \infty ,$ the length of the core 
geodesic $\gamma $ tends to 0. Any sequence of such metrics thus 
converges to the complete (rank 2) hyperbolic cusp
\begin{equation} \label{e4.10}
g_{C} = dr^{2} + e^{2r}g_{o}, 
\end{equation}
on ${\Bbb R} \times T^{2}.$ This process is the formation of a cusp, or 
"opening a cusp", cf. [20], [37].

\begin{remark} \label{r 4.3.}
  {\rm The process described above of opening a cusp may also be 
reversed. Thus, given a complete hyperbolic cusp as in (4.10), the Dehn 
surgery process above closes this cusp by filling in with a hyperbolic 
solid torus, keeping the conformal structure at infinity fixed. As 
discussed above, this can be done in infinitely many non-isometric ways.

 More generally, let $(M^{3}, g_{-1})$ be any complete conformally 
compact hyperbolic 3-manifold with cusps, so that the 
$\varepsilon$-thick part of $M^{3}$ is conformally compact while the 
$\varepsilon$-thin part consists of a finite number of cusps (4.10), 
cf. [37, Ch.5]. Then the Jorgensen-Thurston theory implies that one 
can close the cusps by hyperbolic manifolds, at least for all 
sufficiently short core geodesics, and with at most a small 
perturbation of the structure of conformal infinity.

  In contrast to the situation with solid tori, these manifolds 
obtained by performing Dehn surgery on the cusps of $(M^{3}, g_{-1})$ 
are generally not diffeomorphic. For a fixed diffeomorphism type, 
typically only finitely many such hyperbolic manifolds have a fixed 
conformal infinity.}
\end{remark}

{\bf \S 4.4.}
 The construction for hyperbolic 3-manifolds above is special to 
dimension 3, and cannot be carried out for hyperbolic manifolds in 
dimensions $\geq $ 4. However, we show it can be carried out for AH 
Einstein metrics in dimension 4, (or greater).

 Thus, consider the following $T^{2}$ AdS black hole metrics; we first 
discuss these on the universal cover ${\Bbb R}^{2}\times {\Bbb R}^{2},$ 
and then descend to the quotient ${\Bbb R}^{2}\times T^{2}.$ As before, 
let
\begin{equation} \label{e4.11}
g_{m} = g_{m}^{(0)}= V^{-1}dr^{2} + Vd\theta^{2} + r^{2}(ds_{1}^{2} + 
ds_{2}^{2}), 
\end{equation}
where,
\begin{equation} \label{e4.12}
V = V(r) = r^{2}-  \frac{2m}{r}. 
\end{equation}
As previously, we require $r \geq  r_{+} = (2m)^{1/3} > 0$, where $r_{+}$ 
is the (unique) root of the equation $V(r) = 0$, while $s_{1}, s_{2} 
\in  {\Bbb R} .$ The metric is smooth provided $\theta$ runs over the 
parameter interval [0, $\beta ],$ where $\beta  = \beta_{m}$ is given by
\begin{equation} \label{e4.13}
\beta  = \frac{4\pi}{3r_{+}}, 
\end{equation}

 In contrast to the situation with genus $g \neq 1$ black holes, on 
the space ${\Bbb R}^{2}\times {\Bbb R}^{2},$ the metrics $g_{m}$ are in 
fact all isometric; the change of parameters, (i.e. diffeomorphism), 
given by $r = m^{1/3}s, \theta  = m^{1/3}\psi $ and $s_{i} = 
m^{1/3}t_{i}$ gives an isometry between $g_{m}$ and $g_{1}.$ Thus, in 
the following, we set $m =$ 1.

 Now we essentially repeat the construction in \S 4.3 on these metrics. 
Thus, fix an arbitrary flat structure $g_{o}$ on $T^{2},$ and fix an 
arbitrary simple closed geodesic $\sigma $ in $(T^{2}, g_{o}).$ Let $L 
= L(\sigma )$ be the length of $\sigma $ in $(T^{2}, g_{o}).$ Consider 
first the 3-dimensional metric
\begin{equation} \label{e4.14}
g_{m}'  = V^{-1}dr^{2} + Vd\theta^{2} + r^{2}ds_{1}^{2}, 
\end{equation}
on $D^{2}\times {\Bbb R} ,$ for $V$ as in (4.12) with $m =$ 1. Choose 
$R$ so that
$$V(R)^{1/2}\cdot \beta  = L. $$
Thus, at the boundary $\partial (D^{2}(R)\times {\Bbb R} ),$ the metric 
is the flat metric
\begin{equation} \label{e4.15}
V(R)d\theta^{2} + R^{2}ds_{1}^{2} 
\end{equation}
on the cylinder $S^{1}\times {\Bbb R} .$ The group of Euclidean 
isometries acts on this space, and just as before, given $(T^{2}, 
g_{o})$ there is a unique isometric ${\Bbb Z}$-action on $S^{1}\times 
{\Bbb R} $ such that the ${\Bbb Z}$-quotient metric of (4.15) is 
$(T^{2}, g_{o})$ for which the meridian $\theta $ circle bounding the 
disc is taken to $\sigma .$

 This isometric ${\Bbb Z}$-action on the boundary extends to an 
isometric ${\Bbb Z}$-action on the interior $D^{2}(R)\times {\Bbb R} $ 
and the quotient is a solid torus $D^{2}(R)\times S^{1},$ with 
$\partial D^{2}(R) = \sigma .$ The core geodesic, of distance $R$ to 
the boundary, has length of order $O(e^{-R})$. Further and as before, 
the metric extends to a complete metric on $D^{2}\times S^{1}$.

 In the same way as described in (4.8)-(4.9), one may alter this 
construction slightly to produce such complete, conformally compact 
metrics with the conformal infinity $(T^{2}, g_{o})$ prescribed, in 
place of prescribing the geometry at distance $R$.

 Finally, return to the 4-metric (4.11) and choose an arbitrary, but 
fixed, range for the parameter $s_{2},$ so that $s_{2}\in 
[0,\beta_{2}].$

 To sum up, the analysis above proves:

\begin{proposition} \label{p 4.4.}
  Given any flat structure $(T^{2}, g_{o})$ on the torus, and any 
simple closed geodesic $\sigma $ in $(T^{2}, g_{o}),$ there is a 
complete AH Einstein metric $g$ on the 4-manifold ${\Bbb R}^{2}\times 
T^{2},$ whose conformal infinity is the flat product $(T^{2}, 
g_{o})\times S^{1}(\beta_{2}),$ for any given $\beta_{2} >$ 0. These 
metrics on ${\Bbb R}^{2}\times T^{2}$ are all locally isometric, but 
the isometry type of a metric in this family is uniquely determined by 
the data $(T^{2}, g_{o})$, $\beta_{2}$ and $\sigma$.
\end{proposition}

{\endproof}

 Hence, one has an infinite family of AH Einstein metrics with a given 
conformal infinity. If $\sigma_{i}$ is a sequence of geodesics with 
$L(\sigma_{i}) \rightarrow  \infty ,$ the corresponding metrics $g_{i}$ 
converge to the complete hyperbolic cusp metric   
\begin{equation} \label{e4.16}
g_{C} = dr^{2} + e^{2r}g_{T^{3}}, 
\end{equation}
on ${\Bbb R} \times T^{3},$ where $(T^{3}, g_{T^{3}} = (T^{2}, 
g_{o})\times S^{1}(\beta_{2}).$ Here the convergence is based at points 
$x_{i}$ for which $t_{i}(x_{i}) = 1$ for instance, where $t_{i}$ is the 
geodesic defining function. One thus sees that the regions where the 
metrics $g_{i}$ differ a definite amount from a hyperbolic metric are 
being pushed further and further down the cusp; (this corresponds to 
letting $m \rightarrow $ 0). We also point out that a brief computation 
shows that the $g_{(3)}$ term for any of these metrics satisfies 
$g_{(3)} = 0$; compare with Theorem 2.4.

  These examples illustrate that one may, at least in certain 
situations, open cusps in the class of AH Einstein metrics. Similar but 
more general examples may be obtained by performing Dehn surgery on a 
closed geodesic in $T^{3}$ in place of $T^{2}\times S^{1}$.

\begin{remark} \label{r 4.5.}
  {\rm In analogy to Remark 4.3, it is an interesting open question 
whether this process can be reversed in general. Thus, given a complete 
hyperbolic 4-manifold $M$, with smooth conformal infinity $\partial M,$ 
and with a finite number of cusps, does there exist a sequence of AH 
Einstein manifolds $(M_{i}, g_{i})$, without cusps, such that $(M_{i}, 
g_{i})$ converges to $(M, g)$, and such that the conformal infinity is 
either fixed, or converges to that of $(M, g)$? Again in analogy to 
Remark 4.3, it is to be expected that this requires ${M_{i}}$ to range 
over an infinite collection of topological types in general.} 
\end{remark}

\begin{remark} \label{r 4.6.}
  {\rm All of the discussion in \S 4.1, \S 4.2 and \S 4.4 above 
generalizes in a straightforward way to dimensions $n > $ 4. Thus, one 
replaces the surfaces $\Sigma_{g}, g \geq $ 0, of constant curvature 
$\pm 1,$ 0, by $(n-2)$-dimensional compact Einstein manifolds 
$\Sigma^{n-2}$ of Ricci curvature $\pm (n-3),$ 0. The function $V$ 
becomes $V= c+r^{2}-\frac{2m}{r^{n-3}},$ with $c = \pm 1,$ 0, as 
before.}
\end{remark}

 For the purposes of the next section, for fixed boundary data $(T^{2}, 
g_{o}, \beta_{2}),$ consider the behavior of the geodesic 
compactifications $\bar g_{i}$ on $T^{2}\times {\Bbb R}^{2},$ for 
$g_{i}$ as above with $L(\sigma_{i}) \rightarrow  \infty .$ First, the 
geodesic compactification $\bar g$ of the ${\Bbb R} \times T^{3}$ 
hyperbolic cusp metric (4.16) has the form
\begin{equation} \label{e4.17}
\bar g = dt^{2} + g_{T^{3}}, 
\end{equation}
i.e. $\bar g$ is the flat product metric on ${\Bbb R}^{+}\times T^{3}.$ 
Here of course $t = 2e^{- r},$ and the boundary $\partial M$ occurs at 
$t =$ 0. Note that the ``compactification'' $\bar g$ is not compact, 
due to the cusp. As $i \rightarrow  \infty ,$ the (true) 
compactifications $\bar g_{i}$ converge to $\bar g,$ uniformly on 
compact subsets, based at a point say on $\partial M = \{0\}\times 
T^{3}.$ In particular, (and this is the main point), we have
\begin{equation} \label{e4.18}
diam_{\bar g_{i}}M \rightarrow  \infty , \ \ {\rm as} \ \ i 
\rightarrow  \infty ,
\end{equation}
on $M = {\Bbb R}^{2}\times T^{2}.$

\section{Cusp Formation and Hyperbolic Manifolds.}
\setcounter{equation}{0}

 Proposition 4.4 shows that one may close a complete hyperbolic cusp 
${\Bbb R} \times T^{3}$ in the class of AH Einstein metrics on the 
4-manifold ${\Bbb R}^{2}\times T^{2}$ with a fixed conformal infinity. 
This implies in particular that the space of AH Einstein metrics on a 
fixed manifold $M$ with a fixed conformal infinity is not, in general, 
compact; there are sequences $(M, g_{i})$ of AH Einstein metrics which 
do not converge to an AH Einstein metric on the same space.

 In this section, we prove a type of converse of this statement, namely 
that under reasonable convergence conditions, one can open cusps for AH 
Einstein metrics only when the resulting limit is a complete hyperbolic 
4-manifold. More generally, divergent sequences of AH Einstein metrics 
with controlled conformal infinity can only limit on complete 
hyperbolic 4-manifolds with at least one cusp. The exact statement is 
given in Theorem 5.3.

\medskip

 If $(M, \bar g)$ is a geodesic compactification of an AH Einstein 
manifold $(M, g)$ with geodesic defining function $t$, define its width 
$Wid_{\bar g}(M)$ by
\begin{equation} \label{e5.1}
Wid_{\bar g}(M) = sup\{t(x):x\in M\}. 
\end{equation}
Thus, $Wid_{\bar g}(M)$ is the length of the longest minimizing $\bar g$ 
geodesic starting at $\partial M$ and orthogonal to $\partial M$. Note that 
$Wid_{\bar g}(M)$ depends on the choice of geodesic compactification, 
i.e. the choice of the boundary metric. Two different choices of the 
boundary metric will give rise to different widths, although they can 
be estimated in terms of each other by the conformal factor relating the 
boundary metrics. Observe that the compactifications $\bar g_{i}$ of the 
AH Einstein metrics $g_{i}$ discussed following Remark 4.6 satisfy
\begin{equation} \label{e5.2}
Wid_{\bar g_{i}}(M) \rightarrow  \infty , \ \ {\rm as} \ \ i 
\rightarrow  \infty , 
\end{equation}
corresponding to (4.18).

 As an introduction to the technique, we first show that cusps, (or new 
ends in general), cannot form when the conformal infinity has positive 
scalar curvature.

\begin{proposition} \label{p 5.1.}
  Let $(M, g)$ be an AH Einstein 4-manifold, with boundary metric 
$\gamma .$ Suppose that there is a component of $\partial M$ on which 
the scalar curvature $s_{\gamma}$ of $\gamma $ satisfies
\begin{equation} \label{e5.3}
s_{\gamma} \geq  s_{o} >  0, 
\end{equation}
for some constant $s_{o}.$ Then $\partial M$ is connected and if $\bar 
g$ is the geodesic compactification associated to $\gamma ,$ then
\begin{equation} \label{e5.4}
Wid_{\bar g}M \leq D = \sqrt {3} \pi /\sqrt {s_{o}}. 
\end{equation}
\end{proposition}

\noindent
{\bf Proof:}
 Let $\partial_{o}M$ be a component of $\partial M$ satisfying (5.3). 
For $t_{1} > $ 0 sufficiently small, let $S_{o}(t_{1}) = \{x\in M: 
dist_{\bar g}(x, \partial_{o}M) = t_{1}\},$ so that $S_{o}(t_{1})$ is 
connected and smooth. We may view $S_{o}(t_{1}) \subset $ $(M, g)$, so 
that the function $r$ as in (1.6) has the value $r_{1} = 
\log(\frac{2}{t_{1}})$ on $S_{o}(t_{1}).$

 We construct now a partial defining function $t_{o},$ i.e. a defining 
function for the component $\partial_{o}M$ in the obvious way. Thus, set
$$t_{o} = 2e^{- r_{o}}, $$
where $r_{o}(x) = sgndist(x, S_{o}(t_{1}))+r_{1},$ and $sgndist$ is the 
signed distance function on $(M, g)$ to $S_{o}(t_{1}),$ i.e. 
$sgndist(x) = \pm dist_{g}(x, S_{o}(t_{1}))$ according to whether $t(x) 
<  t_{1}$ or $t(x) >  t_{1}.$

 Note that if $\partial M = \partial_{o}M,$ then $t = t_{o}$ is a 
(full) defining function for the boundary. Otherwise however, $t \neq  
t_{o}$ and the function $t_{o}$ compactifies only the end of $(M, g)$ 
corresponding to $\partial_{o}M,$ in that $\bar g_{o} = t_{o}^{2}g$ is 
compact only on this boundary component. The other boundary components 
of $(M, \bar g_{o})$ are all of infinite $\bar g_{o}$-distance to 
$\partial_{o}M.$

 In either case, it then suffices to show that the maximal length $L$ 
of a (minimizing) $t_{o}$-geodesic $\sigma (t_{o})$ of $\bar g_{o}$ 
satisfies (5.4). This will imply that $\partial M$ is connected, since 
$\partial M$ disconnected implies $L = \infty$, giving a contradiction.
When $\partial M$ is connected, $\bar g_{o} = \bar g$, $L = Wid_{\bar g}(M)$, 
and so (5.4) also follows.

  By (1.19), we have $\bar s_{o}'  = 6t_{o}^{-1}|\bar D^{2}t_{o}|^{2} \geq 0$ 
along $\sigma ,$ so that $\bar s_{o}(\sigma (t_{o}))$ is monotone increasing 
along $\sigma .$ Further, (1.4)-(1.5) and (1.13) imply that 
$Ric_{\bar g_{o}}(N, N) = \frac{1}{6}\bar s_{o} \geq  \frac{1}{4}s_{\gamma}$, 
where $N = \bar \nabla t_{o}$ is the unit tangent vector to $\sigma$. Hence, 
along $\sigma (t_{o})$, one has 
$$Ric_{\bar g_{o}}(N, N) \geq  \frac{1}{4}s_{o}. $$
Now a standard result (Rauch comparison theorem) in Riemannian 
comparison geometry, cf. [35], implies that $\sigma (t_{o})$ must have a 
focal point at distance $D \leq  \sqrt {3} \pi /\sqrt {s_{o}}$, which 
gives (5.4).

  An alternate, even more elementary argument is as follows. The 
equation (1.19) together with the obvious estimate $|\bar 
D^{2}t_{o}|^{2} \geq \frac{1}{3}(\bar \Delta t_{o})^{2}$ and (1.5) imply 
$\bar s_{o}' \geq \frac{1}{18}t_{o} \bar s_{o}^{2}$. Dividing by 
$\bar s_{o}^{2}$ and integrating gives (5.4) with the slightly weaker 
estimate $D = 6/\sqrt{s_{o}}$.
{\endproof}

\begin{remark} \label{r 5.2.}
  {\rm Proposition 5.1 holds in all dimensions, with the same proof, 
cf. also Remark 1.5. As such, it gives a simple new proof of the 
connectedness result of Witten-Yau [39]. More generally, suppose (5.3) 
is replaced by the weaker condition that $s_{\gamma} \geq $ 0 and 
$\partial M$ is not connected. Then there is an infinite 
$t_{o}$-geodesic $\sigma $ of $\bar g_{o}$ joining $\partial_{o}M$ with 
a distinct boundary component of $M$. The argument of Proposition 5.1 
implies that $\bar s_{o} \equiv $ 0 along $\sigma ,$ and hence, via 
(1.19), $D^{2}t_{o} \equiv $ 0 along $\sigma .$ By (1.5), this means 
that $\bar g_{o}$ is Ricci-flat and has a parallel vector field $\nabla 
t_{o}$ along $\sigma ,$ and so the metric $\bar g_{o}$ has an 
infinitesimal splitting as a product of ${\Bbb R} $ with a Ricci-flat 
metric $\gamma_{o}.$ We will see later in Lemma 5.5 that this 
infinitesimal splitting may be globalized to a full splitting, using 
arguments as in the Cheeger-Gromoll splitting theorem. It follows that 
either $\partial M$ is connected or $(M, g)$ is a complete cusp of the 
form
$$g = dr^{2} + e^{2r}\gamma_{o}, $$
where $(\partial_{o}M, \gamma_{o})$ is a compact Ricci-flat manifold. 
This result has been proved by Cai-Galloway [9] using different 
although related methods. Of course, under the weaker bound $s_{\gamma} 
\geq $ 0, even if $\partial M$ is connected one no longer has the 
effective bound (5.4).}
\end{remark}

 We now begin the analysis of the formation of cusps. More generally, 
we study the behavior of sequences $\{g_{i}\}$ of AH Einstein metrics 
on 4-manifolds which have controlled conformal infinities, but which 
diverge in the sense that $\{g_{i}\}$ does not converge to an AH 
Einstein metric on the same manifold.

 Thus, let $(M_{i}, g_{i})$ be a sequence of AH Einstein 4-manifolds, 
with a fixed boundary $\partial M_{i} = \partial M.$ Suppose that 
the conformal infinities $[\gamma_{i}]$ of $g_{i}$ are $C^{m+1, \alpha}$, 
$m \geq 2$, and converge to a limit, so that there are representative metrics 
$\gamma_{i}\in [\gamma_{i}]$ such that $\gamma_{i} \rightarrow \gamma$ 
in $C^{m+1,\alpha}(\partial M)$. Let $\bar g_{i}$ be the associated 
$C^{m, \alpha}$ geodesic compactifications, with $t_{i}$ the associated 
$C^{m+1, \alpha}$ geodesic defining functions. We assume the following:

{\bf Convergence Condition.}
The compactifications $(M_{i}, \bar g_{i})$ converge in the 
$C^{m,\alpha}$ topology, for some $m \geq  2$, and uniformly on {\it compact 
subsets}, to a limit metric $(N, \bar g),$ with boundary metric $\gamma$. 

\smallskip

  This convergence condition should be understood in light of Proposition 
2.7. In particular, the metrics $\bar g_{i}$ and $\bar g$ are smoothed 
near their cutloci to obtain $C^{m, \alpha}$ convergence across the cutlocus. 
It turns out that this convergence condition is not a strong assumption at 
all, but this will only be completely clear in the sequel paper [4].

\medskip

 Since they are distance functions to $\partial M,$ the defining 
functions $t_{i}$ then also converge to the limit $C^{m+1, \alpha}$ geodesic 
defining function $t$ for $\bar g.$ Given base points $x_{i}\in 
t_{i}^{-1}(t_{o}) \subset  M_{i},$ for some fixed $t_{o}$ with 0 $<  
t_{o} <  Wid_{\bar g_{i}}M,$ it follows that the AH Einstein manifolds 
$(M_{i}, g_{i}, x_{i})$ converge, uniformly on compact subsets, to a 
limit complete Einstein manifold $(N, g, x_{\infty})$, $x_{\infty}= \lim 
x_{i},$ with ``compactification'' $\bar g = t^{2}\cdot  g.$ The 
convergence is in the pointed Gromov-Hausdorff topology based at 
$x_{i}$, cf. [21, Ch.3], and also in the $C^{\infty}$ topology, since 
$C^2$ convergence of Einstein metrics implies $C^{\infty}$ convergence, 
by elliptic regularity. 

\medskip

 Now if the width $Wid_{\bar g_{i}}(M_{i})$ of the manifolds $(M_{i}, 
\bar g_{i})$ is uniformly bounded above, it follows by a standard 
application of the Cheeger-Gromov compactness theorem, cf. [2], [12], 
that, in a subsequence, the manifolds $M_{i}$ are all diffeomorphic to 
a fixed manifold $M$, $M = N$, and $\bar g_{i} \rightarrow \bar g$ in 
the $C^{m,\alpha}$ topology on $M$. (The curvature, volume and diameter 
of $(M_{i}, \bar g_{i})$ are all uniformly bounded). Hence, in a 
subsequence, $g_{i}$ is a sequence of AH Einstein metrics on $M$, 
converging to a limit AH Einstein metric $g$ on $M$, for which the 
boundary metrics $\gamma_{i} \rightarrow  \gamma$; compare again with 
Proposition 2.7. In other words, the sequence $(M_{i}, g_{i})$ is not 
divergent in this situation.

 On the other hand, if 
$$Wid_{\bar g_{i}}(M_{i}) \rightarrow  \infty , $$
then any limit complete Einstein manifold $(N, g)$, (again in a 
subsequence), has a non-empty collection of "new" ends, whose boundary 
$\partial_{\infty}N$ is at infinite $\bar g$-distance to 
$\partial N = \partial M.$ In particular, although for any fixed 
$T <  \infty ,$ the domains $U_{i}(T) = t_{i}^{-1}[0,T] \subset  M_{i}$ 
are diffeomorphic to $U (T) = t^{-1}[0,T] \subset  N$, (for $T$ a regular 
value and $i$ sufficiently large), the full manifold $N$ is not diffeomorphic 
to any $M_{i}.$ The discussion concerning and following Proposition 4.4 
exhibits examples where the infinite end of $N$ is a cusp, although this 
of course does not follow automatically in general.

\medskip

  To state the main result on the structure of $(N, g)$ below, we need 
the following two definitions. First, let $\Omega = \Omega(1) = 
t^{-1}[1, \infty) \subset N$ and let $E \subset \Omega$ be any end of 
$\Omega$; (note that $E$ is distinct from an end of $N$ corresponding 
to a boundary component of $\partial M$). Let $S_{E}(t) = S(t) \cap E$, 
where $S(t)$ is the $t$-level set of the geodesic defining function 
$t$ and define
\begin{equation} \label{e5.5}
T_{o}(E) = sup\{t: inf_{S_{E}(t)} \bar s < 0\}.
\end{equation}
If $\bar s \geq 0$ in $E$, set $T_{o}(E) = 0$. Recall again from (1.19) 
that $\bar s$ is non-decreasing along $t$-geodesics in $(\Omega, \bar 
g)$.

  Next, define an end $E \subset \Omega$ as above to be {\it weakly 
hyperbolic} if
\begin{equation} \label{e5.6}
|K + 1|(x) \rightarrow 0, \ \ {\rm as} \ \ t(x) \rightarrow \infty \ \ 
{\rm in} \ \ E,
\end{equation}
where $K$ denotes the sectional curvature of $(E, g)$ at any plane in 
$T_{x}E$.

 Recall also the definition of a conformally compact hyperbolic 
manifold with cusps, as in Remark 4.3, (but in dimension 4 instead of 3). 
We then have the following partial characterization of the limits $(N, g)$ 
obtained above; this result may be considered as a converse to the results 
of \S 4.4, c.f. Proposition 4.4. 

\begin{theorem} \label{t 5.3.}
  Let $(M_{i}, g_{i})$ be a sequence of AH Einstein 4-manifolds, 
$\partial M_{i} = \partial M$, which satisfy the convergence condition. 
Suppose that the Euler characteristics $\chi(M_{i})$ satisfy $\chi 
(M_{i}) \leq  \chi_{o},$ for some $\chi_{o} <  \infty$, that
\begin{equation} \label{e5.7}
Wid_{\bar g_{i}}(M_{i}) \rightarrow  \infty , 
\end{equation}
and that either one of the following two conditions hold:

 (i). There is an end $E \subset \Omega \subset N$ such that $T_{o}(E) 
< \infty$.

 (ii). $(\Omega, g)$ has a weakly hyperbolic end.

 Then the limit $(N, g)$ is a complete conformally compact hyperbolic 
manifold with cusps, with conformal infinity $[\gamma ]$ on $\partial 
M.$ In particular, $(\partial M, [\gamma ])$ is a conformally flat 
3-manifold.
\end{theorem}

 Understandably, the proof is rather long and so is broken into several 
steps. 

\medskip

{\bf Step I.}
First, one needs to control the global size of the AH Einstein manifolds 
$(M_{i}, g_{i})$ away from $\partial M.$

\begin{lemma} \label{l 5.4.}
  Under the assumptions of Theorem 5.3, let $\Omega_{i} = \{x\in 
(M_{i}, \bar g_{i}): t_{i}(x) \geq  1\} = \{x\in (M_{i}, g_{i}): 
r_{i}(x) \leq  \log 2\}$, where $t_{i}$ is the geodesic defining 
function and $r_{i}$ is as in (1.6). Then there is a constant $V^{o} <  
\infty $ such that, for all $i$,
\begin{equation} \label{e5.8}
vol_{g_{i}}\Omega_{i} \leq  V^{o}. 
\end{equation}
\end{lemma}

\noindent
{\bf Proof:}
 For any $(M_{i}, g_{i}),$ the geodesic 'spheres' $S(t) = 
S_{g_{i}}(t),$ i.e. the $t$ level sets of the functions $t_{i},$ have 
the asymptotic expansion (3.5):
\begin{equation} \label{e5.9}
vol_{g_{i}}S(t) = v_{(0)}t^{-3} + v_{(2)}t^{-1}+ o(t). 
\end{equation}
Now by the convergence condition, the geometry of $(M_{i}, \bar g_{i})$ 
between $t =$ 0 and $t =$ 1 is uniformly controlled in $C^{m,\alpha},$ 
and so converges smoothly to that of the limit $(N, \bar g)$ in this 
region. So do the defining functions $t_{i} \rightarrow  t$, and the 
coefficients $v_{(0)}, v_{(2)}.$ Thus, the expansion (5.9) is uniform 
on $S(t)$, in that the lower order term $o(t)$ is small for $t$ small, 
independent of $i$. Hence, for $t_{o}$ small but fixed, by integrating 
(5.9) over the region $t \geq  t_{o},$ we obtain, for $i$ sufficiently 
large,
$$vol_{g_{i}}B(t_{o}) \leq  \frac{1}{3}v_{(0)}t_{o}^{-3} + 
v_{(2)}t_{o}^{-1} + V + 1, $$
where $V = V_{i}$ is the renormalized volume of $(M_{i}, g_{i}),$ cf. 
(3.6), and $B(t_{o}) = t_{i}^{-1}([t_{o}, \infty )) \subset  (M_{i}, 
g_{i}).$ 

 Now by (3.7), the upper bound on $\chi (M_{i})$ gives a uniform upper 
bound on $V$. This gives a uniform upper bound on $vol_{g_{i}}B(t_{o})$ 
and hence (5.8) follows.

{\endproof}

 Summarizing, we have the following description of the structure of the 
limit $(N, g, x_{\infty})$ of $(M_{i}, g_{i}, x_{i}).$ In the collar 
neighborhood $U_{i} = M_{i}\setminus \Omega_{i}$ where $t_{i} \leq $ 1, 
the convergence condition implies that the compactifications $\bar 
g_{i}$ converge smoothly to the compactification $\bar g = t^{2}\cdot  
g$ of the limit. In particular, each $U_{i}$ is diffeomorphic to a 
collar neighborhood $U$ of $\partial M$. By Lemma 5.4, the 
complementary domains $\Omega_{i}$ have uniformly bounded volume, and 
hence the limit region $\Omega  \subset  N$ also has finite volume. 
Further, by (5.7), $Wid_{\bar g}\Omega  = \infty ,$ so that $\Omega  
\subset  N$ has "new" ends, formed from the limiting behavior of 
$(M_{i}, g_{i}).$

  Each end $E$ of $\Omega$ is 'cusp-like' in that it has finite volume, 
and so $volB_{x}(1) \rightarrow 0$, as $x \rightarrow \infty$ in $E$. 
The proof is now split into two cases, according to the hypotheses (i) 
or (ii).

\medskip

 {\bf Step II.}
The following result proves Theorem 5.3 in case (i) holds.
\begin{lemma} \label{l 5.5}
  Suppose that, for some end $E \subset  \Omega ,$
\begin{equation} \label{e5.10}
T_{o}(E) <  \infty . 
\end{equation}
Then $(N, g)$ is a complete conformally compact hyperbolic manifold 
with at least one cusp.
\end{lemma}

\noindent
{\bf Proof:}
Given (5.10), the monotonicity of $\bar s$ implies that there is a 
subend $E'  \subset  E$ on which $\bar s \geq $ 0 and hence there is a 
$t_{o}$ such that $\bar s \geq 0$ on $\Omega_{t_{o}} = \{x \in E: t(x) 
\geq t_{o}\}$. By (1.5), this means that $\bar H = \bar \Delta t \leq 
0$ on $\Omega_{t_{o}}$, where $\bar H$ is the mean curvature of the 
level set $S(t)$, i.e. $S_{E}(t)$, in the direction $\bar \nabla t$. 
Since $Wid_{\bar g}E = \infty$, the discussion in Remark 5.2 shows 
there is an infinitesimal splitting of $(\Omega_{t_{o}}, \bar g)$ along 
a $t$-geodesic ray $\sigma$ in $\Omega_{t_{o}}$.

  To globalize this splitting, consider the domain $\Omega_{t_{o}}$ 
with respect to the Einstein metric $g$. By standard formulas for 
conformal change, cf. also [2,(1.18)], one has
$$H = 3 - t\bar H,$$
where $H$ is the mean curvature of the Lipschitz hypersurface $S(t)$ 
w.r.t. the outward normal $\nabla r$, for $r$ and $t$ related as in 
(1.6). Since $\bar H \leq 0$ on $S(t)$, $t \geq t_{o}$, the mean 
curvature of $(\partial \Omega_{t_{o}}, g)$ satisfies
$$H \geq 3.$$
As in the proof of the Cheeger-Gromoll splitting theorem, cf. 
[6,Ch.6G], this estimate also holds in the sense of distributions or 
support functions at the cut points of $t$ on $S(t_{o})$ where 
$S(t_{o})$ is not smooth. Since we also have $Ric_{g} = -3g$, the 
modification of the Cheeger-Gromoll splitting theorem by Kasue [26], 
cf. also [13], implies that $(\Omega_{t_{o}}, g)$ splits globally as a 
warped product, i.e. as a hyperbolic cusp metric
\begin{equation} \label{e5.11}
g_{C} = dr^{2} + e^{2r}g_{T^{3}}, 
\end{equation}
on ${\Bbb R}\times T^{3}$, where $g_{T^{3}}$ is a flat metric on the 
3-torus $T^{3}$. It follows that the full complete manifold $(N, g)$ is 
hyperbolic, since Einstein metrics are analytic.

{\endproof}

  {\bf Step III.} In this step, we prove Theorem 5.3 in case (ii) 
holds, so that there is an end $E$ which is weakly hyperbolic as in (5.6). 
The next result specifies the geometry of such an end more precisely, 
using Lemma 5.4. 
\begin{lemma} \label{l 5.6}
A weakly hyperbolic end $(E, g)$ of $(N, g)$ is topologically 
${\Bbb R}^{+} \times T^{3}$ and the metric asymptotically approaches 
a hyperbolic cusp metric $g_{C}$, as in (5.11), uniformly on compact sets as 
$t \rightarrow \infty$. More precisely, for any $\varepsilon > 0$ and 
$T < \infty$, there is a $T_{o} = T_{o}(\varepsilon, T)$ such that if 
$t(y) \geq T_{o}$, then the geodesic annulus 
$$A_{y}(T) = \{x \in E: t(x) \in (T^{-1}t(y), Tt(y))\},$$ 
diffeomorphic to $I \times T^{3}$, is $\varepsilon$-collapsed, in that 
$diam_{g}T_{s}^{3} < \varepsilon$, where $T_{s}^{3} = t^{-1}(s)$. 

  Moreover, there exist finite covering spaces $\bar A_{y}(T)$ of $A_{y}(T)$, 
unwrapping the collapse of the $T^{3}$ factors, such that the metric 
$g$ is of the form 
\begin{equation} \label{e5.12}
g|_{\bar A_{y}(T)} = g_{C} + \kappa_{y},
\end{equation}
where the perturbation $\kappa_{y}$ satisfies $||\kappa_{y}|| < \varepsilon$ 
in the $C^{k}$ topology on $\bar A_{y}(T)$, for any given $k < \infty$.
\end{lemma}

\noindent
{\bf Proof:}
The weakly hyperbolic end $(E, g)$ has uniformly bounded curvature, with 
curvature approaching $-1$ as $t \rightarrow \infty$. Further, 
Lemma 5.4 implies that $(E, g)$ has finite volume, so that the volumes of 
unit balls $B_{y}(1)$ tends uniformly to 0 as $t(y) \rightarrow \infty$. 
This means that the manifolds $(E, g, y)$, based at points $y$, are 
collapsing with bounded curvature as $y$ tends to $\infty$ in $E$. Hence, 
the annuli $A_{y}(T)$ have an F-structure, cf. [12], formed essentially by 
the collection of short geodesic loops in $A_{y}(T)$, for $t(y)$ large. 

  When the curvature is highly pinched about $-1$, the structure of such 
collapse is described by the Margulis Lemma, cf. [21], [37]. Thus, as in the 
statement of the Lemma, there are (in fact abelian) covering spaces of the 
annuli $A_{y}(T)$ unwrapping the collapse; the choice of such covering spaces 
is not unique, but they may be chosen so that the injectivity radius and 
diameter of the fibers of the F-structure are on the order of 1. (The 
degree of the covering of course depends on the degree $\varepsilon$ of 
the collapse). In such covering spaces $\bar A_{y}(T)$, the curvature of the 
metric is uniformly close to $-1$, while the diameter and volume of this 
region is uniformly bounded, away from 0 and $\infty$, for $t(y)$ 
sufficiently large. The Cheeger-Gromov compactness theorem then implies 
that the metric $g$ is uniformly close to a hyperbolic metric on 
$\bar A_{y}(T)$. Further, the covering transformations are uniformly close 
to hyperbolic isometries. In any limit as $t(y) \rightarrow \infty$, the 
covering group is hence ${\Bbb Z}^{3}$, the orbits of the F-structure are 
flat $3$-tori, and the limit metric is the hyperbolic cusp metric $g_{C}$ in 
(5.11). Thus, if $t(y)$ is sufficiently large, the metric $g$ on 
$\bar A_{y}(T)$ is uniformly close to a hyperbolic cusp metric. Since the 
metric $g$ is Einstein, the metrics $\bar A_{y}(T)$ are close to $g_{C}$ 
in the $C^{\infty}$ topology.

  Finally, since all geodesic annuli $A_{y}(T)$ are topologically 
$I \times T^{3}$, it follows easily that the end $E$ is topologically 
${\Bbb R}^{+} \times T^{3}$.
{\endproof}

  Lemma 5.6 describes the structure of the end $E$ in the region where 
$t >> 1$. Note that one may let $T \rightarrow \infty$, (sufficiently slowly), 
as $T_{o} \rightarrow \infty$ in Lemma 5.6. In particular, if $y_{k}$ is any 
divergent sequence in $E$, i.e. $t(y_{k}) \rightarrow \infty$ in $(E, g)$, 
then the based sequence $(E, g, y_{k})$ 
has subsequences converging, after unwrapping the collapse as above, to 
the complete hyperbolic cusp metric $g_{C}$. The limit parameter 
$r = r_{\infty}$ in (5.11) is then given by
$$r_{\infty} = \lim_{k\rightarrow\infty}(r -  r(y_{k})),$$ 
with $r = \log(\frac{2}{t})$ as in (1.6). Thus, $r_{\infty}(y) = 0$, where 
$y = \lim y_{k}$ is the limit of the base points $y_{k}$. (This is of course 
analogous to the classical construction of Busemann functions).

 The asymptotic behavior $t >> 1$ of the 'compactification' 
$\bar g = t^{2}\cdot  g$ of $(E, g)$ has a similar description. 
Thus let $t_{k}$ be the geodesic defining function associated with the 
function $r -  r(y_{k}),$ so that
\begin{equation} \label{e5.13}
t_{k} = \frac{1}{t(y_{k})}\cdot t.
\end{equation} 
Thus $t_{k}$ renormalizes $t$ at $y_{k}$, in that $t_{k}$ is a geodesic 
defining function with $t_{k}(y_{k}) = 1$.
The metrics $\bar g_{k} = t_{k}^{2}\cdot  g = (t(y_{k}))^{-2}\cdot \bar g$, 
when based at $y_{k}$, and unwrapped by passing to covers of $T^{3}$ 
as above, have a subsequence converging to the flat product metric
\begin{equation} \label{e5.14}
g_{F} = dt^{2} + g_{T^{3}}, 
\end{equation}
on $F = {\Bbb R}^{+}\times T^{3}$; this follows for example from the formulas 
(1.3)-(1.5) and (1.19), compare with (4.17). Here, the limit parameter 
$t$ is given by $t \equiv t_{\infty} = \lim_{k\rightarrow\infty}t_{k},$ 
associated to $r_{\infty}$ as above. Of course for $y = \lim y_{k}$ as 
above, $t(y) = t_{\infty}(y) = 1$, so that $y \in F$ has distance 1 to 
$\partial F = \{0\} \times T^{3}$.

 This discussion holds for {\it  any}  divergent sequence $\{y_{k}\}$ 
in $E$. Note however that we do not assert that the flat structure on 
$T^{3}$ is independent of the sequence $\{y_{k}\}.$ Apriori it is possible 
that different sequences may give rise to flat limits (5.14) with distinct 
flat structures on $T^{3}$, although if $y_{k}$ and $y_{k}' $ are distinct 
sequences with $t(y_{k})/t(y_{k}')$ bounded away from 0 and infinity, 
then the limit metrics are the same, (i.e. isometric). This possibility of 
the non-uniqueness of the 'tangent cones at infinity', does not play any 
role however in the remainder of the proof. 

  It is worth emphasizing again that, for any divergent sequence $\{y_{k}\}$, 
$t(y_{k}) \rightarrow \infty$ in $(E, g)$, the 
sequence of metrics $\bar g_{k}$ as $k \rightarrow  \infty $ describes 
the normalized asymptotic behavior in regions about $y_{k}$ of the 
{\it  fixed}  metric $(E, \bar g)$, in that the metrics 
$\{\bar g_{k}\}$ are just rescalings and unwrappings of $\bar g$ based 
at $y_{k}$.

\medskip

An end $E \subset \Omega$ having the structure described in Lemma 5.6 will 
be called an {\it asymptotically hyperbolic cusp}. Theorem 5.3 is now an 
immediate consequence of the following rigidity result.

\begin{proposition} \label{p 5.7.}
  Let $(N, g)$ be an AH Einstein 4-manifold, with at least one 
asymptotically hyperbolic cusp $E$. Then $(N, g)$ is hyperbolic.
\end{proposition}

\noindent
{\bf Proof:}
 As described above, for {\it any} divergent sequence $y_{k}$ in $E$, the 
Riemannian manifolds $(E, g, y_{k})$ converge, in a subsequence and 
uniformly on compact subsets, to a complete hyperbolic cusp after unwrapping 
the collapse. The limit parameter $r = r_{\infty}$ is normalized by 
$r(y) = 0$, where $y$ is the limit base point. Thus, $\{y_{k}\}$ determines a 
sequence of Einstein perturbations of the hyperbolic cusp metric (5.11). 
If $(E, g)$ itself is not hyperbolic, then the based metrics $(E, g, y_{k})$ 
are not hyperbolic, so that the sequence of perturbations is non-trivial. 
We will prove that this assumption leads to a contradiction.

 For computation, it is convenient, (although not necessary), to work 
with the compactification $\bar g.$ Thus, as described above, the 
compactifications $\bar g_{k} = t_{k}^{2}\cdot  g$ based at $y_{k}$ 
converge, in a subsequence, to a flat product metric $g_{F}$ (5.14) on 
${\Bbb R}^{+}\times T^{3},$ again after unwrapping the collapse. The 
convergence of $\bar g_{k}$ to $g_{F}$ is smooth and uniform on compact 
subsets of ${\Bbb R}^{+}\times T^{3},$ but is not smooth at the 
boundary $\{0\}\times T^{3}.$ 

 Now view the metrics $\bar g_{k}$ as perturbations of the limit flat 
metric $g_{F}.$ Note that the metrics $\bar g_{k}$ are all Bach-flat, 
i.e. satisfy the Bach equation (2.5). If any $\bar g_{k}$ is flat on 
some open set $U \subset  E$ containing some $y_{k}$, then $g$ is locally 
conformally flat in $U$. Since $g$, being Einstein, is analytic, it is then 
everywhere locally conformally flat and hence $(N, g)$ is hyperbolic, i.e. 
the result follows in this case. Thus, we may and do assume that $\bar g_{k}$ 
is not flat on any open set, for all $k$.

 To understand the behavior of $\bar g_{k}$ near the flat limit 
$g_{F},$ consider the linearization. Thus write
\begin{equation} \label{e5.15}
\bar g_{k} = g_{F} + s_{k}h_{k}, 
\end{equation}
where $s_{k} \rightarrow 0$ and $h_{k}$ is a sequence of symmetric 
bilinear forms with $s_{k}h_{k} \rightarrow 0$ smoothly on compact subsets. 
As above, it is understood here and below that the metrics $\bar g_{k}$ are 
lifted to covering spaces unwrapping the collapse. The parameter $s_{k}$ 
is chosen measure the local size of the curvature at the base point $y_{k}$ 
in that
\begin{equation} \label{e5.16}
s_{k} = (\int_{B_{y_{k}}(\frac{1}{2})}|R_{\bar g_{k}}|^{2}dV)^{1/2}. 
\end{equation}
Since $\bar g_{k}$ is not flat anywhere, $s_{k} > $ 0, for all $k$. The 
convergence $\bar g_{k} \rightarrow  g_{F}$, (in a subsequence), is smooth, 
and so the forms $h_{k}$ are locally bounded, away from $\{0\}\times T^{3}$, 
and converge smoothly to a limit symmetric bilinear form $h$ on 
${\Bbb R}^{+}\times T^{3}$, with $||h|| \sim 1$ at the base point 
$y = \lim y_{k}$. 
Further, since the convergence of $\bar g_{k}$ to $g_{F}$ requires unwrapping 
to larger and larger covers, the limit form $h$ is invariant under the $T^{3}$ 
action on ${\Bbb R}^{+}\times T^{3}.$

 The limit $h$ is not uniquely defined, since one may alter the 
convergence $\bar g_{k} \rightarrow  g_{F}$ by diffeomorphisms 
converging to the identity; this corresponds to changing $h$ to 
$h+\delta^{*}X,$ for some vector field $X$. To normalize, $h$ may be 
chosen so that
\begin{equation} \label{e5.17}
\beta_{g_{F}}(h) = 0, 
\end{equation}
where $\beta_{g_{F}}$ is the Bianchi operator of $g_{F}, 
\beta_{g_{F}}(h) = \delta h + \frac{1}{2}dtrh,$ where the divergence 
and trace are w.r.t $g_{F},$ cf. also [7]. 

  Now the form $h$ is a solution of the linearized Bach equations at the 
flat metric $g_{F}$ and the deviation of $\bar g_{k}$ from $g_{F}$ is 
measured, to first order, by the size of the linearization $h$, in that 
\begin{equation} \label{e5.18}
\bar g_{k} = g_{F} + s_{k}h + o(s_{k}),
\end{equation}
where $o(s_{k}) << s_{k}$ on any given compact subset of 
${\Bbb R}^{+} \times T^{3}$. In particular, the curvature $R_{\bar g_{k}}$ 
on the annuli $A_{y_{k}}(T)$ satisfies
\begin{equation} \label{e5.19}
|R_{\bar g_{k}}| \sim s_{k} |\partial^{2}h + Q_{1}(h)|,
\end{equation}
where $Q_{1}(h)$ involves only $h$ and its first derivative. Note that 
$|R_{\bar g_{k}}| \sim s_{k}$ on the $L^{2}$-average in 
$B_{y_{k}}(\frac{1}{2})$, by (5.16). The following Lemma gives the 
structure of any such linearization $h$ which arises from an Einstein 
perturbation of a hyperbolic cusp metric, as above.
\begin{lemma} \label{l 5.8}
Any $T^{3}$ invariant symmetric bilinear form $h$ on 
${\Bbb R}^{+} \times T^{3}$ constructed as above and satisfying (5.17) is 
given by
\begin{equation} \label{e5.20}
h = C^{(0)}+ C^{(1)}t + C^{(2)}t^{2}+ C^{(3)}t^{3}+ C^{(4)}t^{4}, 
\end{equation}
where the coefficients $C^{(i)}$ are constant, i.e. parallel forms, on 
${\Bbb R}^{+} \times T^{3}$, and $t = t_{\infty}$ is the parameter on 
${\Bbb R}^{+}$, as in (5.14).
\end{lemma}

\noindent
{\bf Proof:}
Since the limit is flat, it is easily seen from (2.5) that the linearized 
Bach equation is
\begin{equation} \label{e5.21}
2D^{*}D(Ric' (h)) = -\frac{2}{3}D^{2}s'  -\frac{1}{3}(\Delta s' )g_{F}, 
\end{equation}
where $Ric' (h) = \frac{d}{ds}Ric(g_{F}+sh)$ is the linearization of 
the Ricci curvature at the flat metric $g_{F}$ and similarly $s'  = s' 
(h)$ is the linearization of the scalar curvature, in the direction 
$h$. From standard formulas, cf. [6, Ch.1K], the normalization (5.17) 
at the flat metric gives 
\begin{equation} \label{e5.22}
Ric' (h) = \frac{1}{2}D^{*}Dh, \ \ {\rm and} \ \ s' (h) = 
-\frac{1}{2}\Delta trh. 
\end{equation}
Hence, (5.21) becomes
\begin{equation} \label{e5.23}
(D^{*}D)^{2}h = \frac{1}{3}D^{2}(\Delta trh) + \frac{1}{6}(\Delta\Delta 
trh)g_{F}. 
\end{equation}

 The task now is to determine the $T^{3}$ invariant solutions of 
(5.23). To do this, let $e_{i}$ be an orthonormal framing for the flat 
metric $g_{F},$ with $e_{1} = \nabla t,$ and $e_{i}, i =$ 2,3,4 tangent 
to the $T^{3}$ factor and let $\theta_{i}$ be the corresponding 
coframing. Thus
$$h = \sum h_{ij}\theta_{i}\cdot \theta_{j}, $$
where $h_{ij} = h_{ij}(t),$ since $h$ is $T^{3}$ invariant. It is 
straightforward to compute that the Bianchi normalization (5.17) gives 
the equations
\begin{equation} \label{e5.24}
\partial_{t}h(e_{1},e_{1}) = \frac{1}{2}\partial_{t}trh, \ \  
\partial_{t}h(e_{1},e_{i}) = 0, i \geq  2. 
\end{equation}
Next, recall from (1.5) that the scalar curvature of a geodesic 
compactification is given by $s = - 6\frac{\Delta t}{t};$ (as usual we 
drop the overbars). Hence, $s'  = 6t^{-2}(\Delta t)t'  -  
6t^{-1}(\Delta' )(t) -  6t^{-1}\Delta (t')$, where $t'$ is the 
linearization of $t$ in the direction $h$. The first term here vanishes, 
since $\Delta t = 0$ on $g_{F}.$ For the second term, from [6,Ch.1K], 
$(\Delta')(t) = -< D^{2}t, h>  + < dt, \beta (h)>  = 0$, by (5.17) and 
the fact that $D^{2}t = 0$ on $g_{F}$. Thus,  
\begin{equation} \label{e5.25}
s'  = - 6t^{-1}\Delta (t') = - 6t^{-1}\partial_{t}\partial_{t}(t'); 
\end{equation}
here the second equality follows from the fact that $t'$ is only a function of 
$t$, since $h$ is. To compute $\partial_{t}(t')$, let $g_{s} = g_{F} + sh$ and 
let $t_{s}$ be distance functions w.r.t. $g_{s}$ converging to the distance 
function $t = t_{\infty}$ on $(F, g_{F})$. (For example for $s = s_{k}$ as in 
(5.18), $t_{s} = t_{s_{k}} = t_{k}$ is given as in (5.13)). We have 
$t_{s} = t + st' + o(s)$ and $|\nabla_{\bar g_{s}} t_{s}|^{2} = 1$, 
i.e. $\bar g_{s}^{ij}\partial_{i}t_{s}\partial_{j}t_{s} = 1$. Taking the 
derivative w.r.t. $s$ then gives, at $g_{F}$,
$$2\partial_{t}(t') = 2<\nabla t' , \nabla t>  = h(\nabla t, \nabla t) 
= h(e_{1},e_{1}). $$
Combining this with (5.24) and (5.25) results in
\begin{equation} \label{e5.26}
s'  = -\frac{3}{2}t^{-1}\partial_{t}trh, 
\end{equation}
which, combined with (5.22) gives $\partial_{t}\partial_{t}trh = 
3t^{-1}\partial_{t}trh$. Hence, $\partial_{t}trh = c_{o}t^{3}$ and so
\begin{equation} \label{e5.27}
trh = \frac{c_{o}}{4}t^{4} + c_{1}, 
\end{equation}
for some constants $c_{o}, c_{1}.$ Thus (5.23) reduces to the $4^{\rm th}$ 
order equation
$$h_{ij}^{(iv)} = 2c_{o}\delta_{1i}\delta_{1j} + c_{o}\delta_{ij},$$ 
which implies the result.
{\endproof}

  The polynomials of order 1, i.e. the constant and linear forms in $t$ in 
(5.20), give rise to trivial, i.e. flat deformations of $g_{F}$, and so 
don't contribute to the curvature in (5.19) or (5.22). Since 
by construction, i.e. by the choice of $s_{k}$ in (5.16), $h$ is non-trivial, 
$h$ contains polynomials of degree at least 2. Hence
\begin{equation} \label{e5.28}
\nabla^{2}h = \nabla_{T}\nabla_{T}h = 2C^{(2)} + 6C^{(3)}t + 12C^{(4)}t^{2}.
\end{equation}
Observe that in the context of the perturbation $\bar g_{k}$ in (5.15), 
(5.19ff) implies that $C^{(2)}$ is uniformly bounded away from 0 and 
$\infty$. 

  Now return to the geometry of the metric $(E, g)$ or $(E, \bar g)$, (with 
the original defining function $t$, in place of $t = t_{\infty}$ above). 
Observe that the results above hold for {\it any} divergent sequence of base 
points $\{y_{k}\}$ in $(E, \bar g)$, i.e. $t(y_{k}) \rightarrow \infty$. 
This means that for any $y \in E$ with $t(y)$ sufficiently large, 
the rescaled metrics $\bar g_{y} = t_{y}^{2} \cdot g$, 
$t_{y} = t/t(y)$, based at $y$, are always of the form (5.18) on large annuli, 
with $h = h_{y}$ of the form (5.20). As noted following (5.13), observe 
also that $\bar g_{y} = t(y)^{-2} \bar g$, i.e. $\bar g_{y}$ is a constant 
(not conformal) rescaling of $\bar g$.

  Now on the one hand, by the weakly hyperbolic assumption (5.6), all the rescaled 
metrics $\bar g_{y}$ tend to the flat metric as $t(y) \rightarrow \infty$, (i.e. the 
parameter $s_{k}$, now $s_{y}$, in (5.16) tends to 0). On the other hand, (5.28) and (5.19) 
show that, at any given $y$, the curvature of $\bar g_{y}$, although of 
necessity small, at least remains bounded away from 0 on 
$(A_{y}(T), \bar g_{y})$, for $T >> 1$ and  $t_{y}$ large. In other words, 
for any $y' \in (A_{y}(T), \bar g_{y})$ with $t(y') >> t(y)$, one has
$$|R_{\bar g_{y}}|(y') \geq c_{o} |R_{\bar g_{y}}|(y),$$
where $c_{o}$ is a fixed numerical constant. For any such $y'$ since 
$\bar g_{y'} = (t(y)/t(y'))^{2}\bar g_{y}$, and 
$t(y)/t(y') << 1$, it follows that
\begin{equation} \label{e5.29}
|R_{\bar g_{y'}}|(y') = (t(y')/t(y))^{2}|R_{\bar g_{y}}|(y') >> |R_{\bar g_{y}}|(y)|.
\end{equation}
The estimate (5.29) implies, for instance by iteration, that the curvature of 
$\bar g_{y}$ cannot decrease to 0 as $t(y) \rightarrow \infty$ in $E$. This 
is of course a contradiction.

 This contradiction implies that there are no non-trivial $T^{3}$-invariant 
Bach-flat deformations of the flat metric arising in this way, and hence no 
non-trivial deformations of the hyperbolic cusp metric among Einstein 
metrics. As explained at the beginning of the proof, this contradiction 
proves Proposition 5.7, which thus also completes the proof of Theorem 
5.3.

{\endproof}

  It is an interesting open question whether Theorem 5.3 remains valid 
without one of the hypotheses (i) or (ii).

\bibliographystyle{plain}

\begin{thebibliography}{WWW}

\footnotesize


\bibitem[1]{1} R. A. Adams, ``Sobolev Spaces'', Pure and Appl. Math. Series, 
vol. 65, Academic Press, New York, (1975).

\bibitem[2]{2} M. Anderson, Extrema of curvature functionals on the 
space of metrics on 3-manifolds, {\it Calc. Var. and P.D.E.}, {\bf 5}, 
(1997), 199-269.

\bibitem[3]{3} M. Anderson, $L^{2}$ curvature and volume 
renormalization for AHE metrics on 4-manifolds,  {\it Math. Research Lett.}, 
{\bf 8}, (2001), 171-188: math-DG/0011051.

\bibitem[4]{4} M. Anderson, Einstein metrics with prescribed conformal 
infinity on 4-manifolds, (preprint, May 2001): math.DG/0105243.

\bibitem[5]{5} M. Berger, Quelques formules de variation pour une 
structure Riemannienne, {\it Ann. Sci. Ecole Norm. Sup.}, {\bf 3}, (1970), 
285-294.

\bibitem[6]{6}A. Besse, ``Einstein Manifolds'', Ergebnisse Series, vol 
3:10, Springer Verlag, New York, (1987).

\bibitem[7]{7}O. Biquard, ``M\'etriques d'Einstein asymptotiquement 
sym\'etriques'', Ast\'erisque, {\bf 265}, (2000).

\bibitem[8]{8}O. Biquard, Einstein deformations of hyperbolic metrics, 
{\it} in ``Essays on Einstein Manifolds'', Surveys in Diff. Geom. VI, Ed. C. 
LeBrun and M. Wang, Intern. Press, Cambridge, MA, (1999), 235-246.

\bibitem[9]{9}M. Cai and G.J. Galloway, Boundaries of zero scalar 
curvature in the AdS/CFT correspondence, {\it Adv. Theor. Math. Phys.}, 
{\bf 3} (1999), 1769-1783: hep-th/0003046

\bibitem[10]{10}A.P. Calder\'on, Uniqueness in the Cauchy problem for 
partial differential equations, {\it Amer. Jour. Math.}, {\bf 80}, (1958), 
16-36.

\bibitem[11]{11}A.P. Calder\'on, Existence and uniqueness theorems for 
systems of partial differential equations, {\it in} ``Proc. Symp. Fluid 
Dynamics and Appl. Math'', (Univ. of Maryland, 1961), Gordon and Breach, N.Y., 
(1962), 147-195.

\bibitem[12]{12}J. Cheeger and M. Gromov, Collapsing Riemannian 
manifolds while keeping their curvature bounded, I, II, {\it Jour. Diff. 
Geom.}, {\bf 23}, (1986), 309-346, and {\bf 32}, (1990), 269-298.

\bibitem[13]{13}C.B. Croke and B. Kleiner, A warped product splitting 
theorem, {\it Duke Math. Jour.}, {\bf 67}, (1992), 571-574.

\bibitem[14]{14}S. de Haro, K. Skenderis and S.N. Solodukhin, 
Holographic reconstruction of spacetime and renormalization in the 
AdS/CFT correspondence, {\it Comm. Math. Phys.}, {\bf 217}, (2001), 595-622, 
hep-th/0002230.

\bibitem[15]{15} J. Escobar, On the prescribed scalar curvature problem 
on compact manifolds with boundary, {\it Contemp. Math.}, {\bf 268}, Amer. 
Math. Soc., (2000), 137-144.

\bibitem[16]{16}C. Fefferman and C.R. Graham, Conformal invariants, {\it in} 
``\'Elie Cartan et les Mathematiques d'Aujourd'hui'', Ast\'erisque, (1985), 
95-116.

\bibitem[17]{17}D. Gilbarg and N. S. Trudinger, ``Elliptic Partial 
Differential Equations of Second Order'', $2^{\rm nd}$ Edition, Springer 
Verlag, New York, 1983.

\bibitem[18]{18}C.R. Graham, Volume and area normalization for 
conformally compact Einstein metrics, {\it Rend. Circ. Mat. Palermo}, 
Ser II, Suppl., {\bf 63}, 2000, 31-42: math.DG/0009042

\bibitem[19]{19}C.R. Graham and J.M. Lee, Einstein metrics with 
prescribed conformal infinity on the ball, {\it Advances in Math.}, {\bf 87}, 
(1991), 186-225.

\bibitem[20]{20}M. Gromov, Hyperbolic manifolds according to Thurston 
and Jorgensen, {\it Seminaire Bourbaki}, No. 546, (1980).

\bibitem[21]{21}M. Gromov, ``Metric Structures for Riemannian and 
Non-Riemannian Spaces'', Prog. in Math. Series, {\bf 152}, Birkhauser 
Verlag, Boston, (1999).

\bibitem[22]{22}M. Gursky, Compactness of conformal metrics with 
integral bounds on curvature, {\it Duke Math. Jour.}, {\bf 72}, (1993), 
339-367.

\bibitem[23]{23}S. W. Hawking and D. N. Page, Thermodynamics of black 
holes in Anti-de Sitter space, {\it Comm. Math. Phys.}, {\bf 87}, (1983), 
577-588.

\bibitem[24]{24}S. W. Hawking, C.J. Hunter and D. N. Page, Nut charge, 
Anti-de Sitter space and entropy, {\it Phys. Rev. D}, {\bf 59}, (1999), 
044033: hep-th/9809035.

\bibitem[25]{25}M. Henningson and K. Skenderis, The holographic Weyl 
anomaly, {\it Jour. High Energy Phys.}, {\bf 9807}, (1998), 023, 
hep-th/9806087.

\bibitem[26]{26}A. Kasue, Ricci curvature, geodesics and some geometric 
properties of Riemannian manifolds with boundary, {\it J. Math. Soc. Japan}, 
{\bf 35}, (1983), 117-131.

\bibitem[27]{27}S. Kobayashi and K. Nomizu, ``Foundations of Differential 
Geometry'', vol. 1, Interscience, J.W. Wiley, New York, (1963).

\bibitem[28]{28}C. LeBrun, ${\cal H}$-space with a cosmological 
constant, {\it Proc. Royal Soc. London}, Ser. A, {\bf 380}, (1982), 171-185.

\bibitem[29]{29} M. Li, The Yamabe problem with Dirichlet data, {\it C. R. 
Acad. Sci. Paris}, {\bf 320}, S\'erie I, (1995), 709-712.

\bibitem[30]{30}J.L. Lions and E. Magenes, ``Non-Homogeneous Boundary 
Value Problems and Applications'', I, Grundlehren Series, vol. 181, 
Springer Verlag, New York, 1972.

\bibitem[31]{31}J. Maldacena, The large $N$ limit of superconformal 
field theories and supergravity, {\it Adv. Theor. Math. Phys}, {\bf 2}, 
(1998), 231-252, hep-th/9711200

\bibitem[32]{32}C.B. Morrey, Jr. ``Multiple Integrals in the Calculus of 
Variations'', Grundlehren Series, vol. 130, Springer Verlag, Berlin, 1966.

\bibitem[33]{33}R. C. Myers, Stress tensors and Casimir energies in the 
AdS/CFT correspondence, (preprint): hep-th/9903203.

\bibitem[34]{34}L. Nirenberg, ``Lectures on Partial Differential 
Equations'', CMBS-NSF Reg. Conf. Series in Math., {\bf 17}, Amer. Math. 
Soc., (1973).

\bibitem[35]{35}P. Petersen, ``Riemannian Geometry'', Grad. Texts in Math., 
vol. 171, Springer Verlag, New York, (1997).

\bibitem[36]{36}K.T. Smith, Some remarks on a paper of Calderon on 
existence and uniqueness theorems forsystems of partial differential 
equations, {\it Comm. Pure Appl. Math.}, {\bf 18}, (1965), 415-441.

\bibitem[37]{37}W. Thurston, ``The Geometry and Topology of 
Three-Manifolds'', (preprint), Princeton Univ.,(1978)

\bibitem[38]{38}E. Witten, Anti De Sitter space and holography, {\it Adv. 
Theor. Math. Phys.}, {\bf 2}, (1998), 253-291, hep-th/9802150

\bibitem[39]{39}E. Witten and S.-T. Yau, Connectedness of the boundary 
in the AdS/CFT correspondence, {\it Adv. Theor. Math. Phys.}, {\bf 3}, 
(1999), 1635-1655, hep-th/9910245.

\end{thebibliography}

\bigskip

\begin{center}
April, 2001/May, 2002
\end{center}

\medskip

\noindent
\address{Department of Mathematics\\
S.U.N.Y. at Stony Brook\\
Stony Brook, N.Y. 11790-3651}\\
\email{anderson@@math.sunysb.edu}

\end{document}